\documentclass{article}
\usepackage{amsfonts,amsmath,amssymb}
\usepackage{graphics, epsfig}
\usepackage{color}
\usepackage{appendix}
\allowdisplaybreaks

\newtheorem{theorem}{Theorem}[section]
\newtheorem{proposition}{Proposition}[section]
\newtheorem{lemma}{Lemma}[section]
\newtheorem{corollary}{Corollary}[section]
\newtheorem{remark}{Remark}[section]

\numberwithin{equation}{section}
\numberwithin{theorem}{section}
\numberwithin{proposition}{section}
\numberwithin{lemma}{section}
\numberwithin{remark}{section}
\numberwithin{corollary}{section}
\newcommand{\noi}{\noindent}

\newcommand{\dsty}{\displaystyle}
\newcommand{\txty}{\textstyle}

\newcommand{\al}{\alpha}

\newcommand{\be}{\beta}

\newcommand{\gm}{\gamma}
\newcommand{\dl}{\delta}
\newcommand{\Dl}{\Delta}
\newcommand{\lm}{\lambda}
\newcommand{\Lm}{\Lambda}

\newcommand{\varep}{\varepsilon}
\newcommand{\eps}{\epsilon}
\newcommand{\vp}{\varphi}
\newcommand{\sig}{\sigma}

\newcommand{\z}{\zeta}

\newcommand{\df}[1]{\buildrel\mbox{\small def}\over{#1}}

\newcommand{\db}{\prime\prime}

\newcommand{\nn}{\mathbb{N}}

\newcommand{\rr}{\mathbb{R}}
\newcommand{\rn}{\rr^N}

\newcommand{\mcl}[1]{\mathcal{#1}}
\newcommand{\bl}[1]{\mathbf{#1}}

\newcommand{\bbox}{\vrule height.6em width.6em 
depth0em} 

\newcommand{\dvg}{\operatorname{div}}

\newcommand{\essup}{\operatornamewithlimits{ess\,sup}}

\newcommand{\loc}{\operatorname{loc}}

\newcommand{\bsu}{\mathop{\txty{\sum}}\limits}

\newcommand{\pl}{\partial}

\newcommand{\intl}{\int\limits}

\def\Xint#1{\mathchoice
    {\XXint\displaystyle\textstyle{#1}}%
    {\XXint\textstyle\scriptstyle{#1}}%
    {\XXint\scriptstyle\scriptscriptstyle{#1}}%
    {\XXint\scriptscriptstyle\scriptscriptstyle{#1}}%
    \!\int}
\def\XXint#1#2#3{\setbox0=\hbox{$#1{#2#3}{\int}$}
    \vcenter{\hbox{$#2#3$}}\kern-0.5\wd0}
\def\bint{\Xint-}
\def\dashint{\Xint{\raise4pt\hbox to7pt{\hrulefill}}}

\newcommand{\ovl}[3]{\int_{#1}^{#2}\kern-#3pt\raise4pt\hbox to7pt{\hrulefill}\ }

\newcommand{\ovll}[3]{\intl_{#1}^{#2}\kern-#3pt\raise4pt\hbox to7pt{\hrulefill}\ }

\newcommand{\tvl}[2]{\iint_{#1}\kern-#2pt\raise4pt\hbox to7pt{\hrulefill}\ }

\newcommand{\bye}{\end{document}}

\newcommand{\testfunlu}{\frac{u^{-\frac m2}-\ssigbar^{-\frac m2}}{m}}
\newcommand{\testfunls}{\frac{s^{-\frac m2}-\ssigbar^{-\frac m2}}{m}}
\newcommand{\testfunmu}{\frac{\ssigbar^{-\frac m2}-u^{-\frac m2}}{m}}
\newcommand{\testfunms}{\frac{\ssigbar^{-\frac m2}-s^{-\frac m2}}{m}}

\newcommand{\data}{\{N,C_o,C_1\}}
\newcommand{\pto}{(x_o,t_o)}

\newcommand{\tkn}{\tilde{K}_n}
\newcommand{\trn}{\tilde{\rho}_n}

\newcommand{\tvls}[2]{\iint_{#1}\kern-#2pt\raise4pt\hbox to15pt{\hrulefill}\ }

\newcommand{\uqoo}{\bint_{K_\rho(x_o)}u^q(\cdot,t_o)dx} 
\newcommand{\uqoonrm}{\Big(\uqoo\Big)^{\frac1q}}

\newcommand{\ssig}{\mcl{S}_\sig}
\newcommand{\ssigbar}{\mcl{\bar S}_\sig}
\begin{document}
\title{Two Remarks on the Local Behavior of Solutions 
to Logarithmically Singular Diffusion Equations and 
its Porous-Medium Type Approximations}
\author{Emmanuele DiBenedetto\\
Department of Mathematics, Vanderbilt University\\  
1326 Stevenson Center, Nashville TN 37240, USA\\
email: {\tt em.diben@vanderbilt.edu}
\and
Ugo Gianazza\\
Dipartimento di Matematica ``F. Casorati", 
Universit\`a di Pavia\\ 
via Ferrata 1, 27100 Pavia, Italy\\
email: {\tt gianazza@imati.cnr.it}
\and
Naian Liao\\
Department of Mathematics, Vanderbilt University\\  
1326 Stevenson Center, Nashville TN 37240, USA\\
email: {\tt naian.liao@vanderbilt.edu}}
\date{}
\maketitle
\vskip.4truecm
\begin{abstract}
For the logarithmically singular parabolic 
equation \eqref{Eq:1:1} below, we establish a Harnack 
type estimate in the $L^1_{\loc}$ topology,
and we show that the solutions are locally analytic in the 
space variables and differentiable in time. 
The main assumption is that $\ln u$ possesses a 
sufficiently high degree of integrability 
(see \eqref{Eq:1:3} for a precise statement). These 
two properties are known for solutions of 
singular porous medium type equations ($0<m<1$), 
which formally approximate the logarithmically 
singular equation \eqref{Eq:1:1} below. However, the corresponding 
estimates deteriorate as $m\to0$. It is 
shown that these estimates become stable and 
carry to the limit as $m\to0$, provided the indicated 
sufficiently high order of integrability is in force. 
The latter then appears as the discriminating 
assumption between solutions of parabolic 
equations with power-like singularities and 
logarithmic singularities to insure such 
solutions to be regular.  
\end{abstract}
\section{Main Results}\label{S:1}
We continue here the investigations initiated in \cite{DBGL1,DBGL2}, 
on the local behavior of non-negative solutions to logarithmically 
singular parabolic equations of the type
\begin{equation}\label{Eq:1:1}
\begin{aligned}
& u\in C_{\loc}\big(0,T; L^2_{\loc}(E)\big),\quad 
\ln u\in L^2_{\loc}\big(0,T;W^{1,2}_{\loc}(E)\big);\\
&u_t-\Delta \ln u=0 \qquad\text{ weakly in }\> E_T=E\times(0,T]
\end{aligned}
\end{equation}
where $E$ is an open set in $\rn$ and $T>0$.  It is 
assumed throughout that 
\begin{equation}\label{Eq:1:2}
u\in L^r_{\loc}(E_T)\quad\text{ for some } 
\> r>\max\big\{1;{\txty\frac N2}\big\}
\end{equation}
and that 
\begin{equation}\label{Eq:1:3}
\ln u\in L^\infty_{\loc}\big(0,T;L^p_{\loc}(E)\big)
\quad\text{ for some } 
\> p\ge1.  
\end{equation}
The modulus of ellipticity of the principal part is 
$u^{-1}$. Therefore the equation is degenerate as 
$u\to\infty$ and singular as $u\to0$. It was shown 
in \cite{DBGL1} that (\ref{Eq:1:2}) implies that 
$u$ is locally bounded in $E_T$, and hence the 
equation is not degenerate. Likewise if  (\ref{Eq:1:3}) 
holds for some $p>N+2$, then the solution $u$ is locally 
bounded below, and hence the equation 
is not singular. As a consequence $u$ is locally, 
a classical solution to (\ref{Eq:1:1}). This was 
realized by establishing a local upper and lower 
bound on $u$, via a pointwise Harnack-type estimate. 

The main results of this note are that if $u$ is 
a locally bounded, weak solutions to (\ref{Eq:1:1}), 
then:  
\begin{description}
\item{\bf i.\ } If $\ln u$  satisfies (\ref{Eq:1:3}) for 
some $p\ge2$, then $u$ satisfies a local 
Harnack inequality in the $L^1_{\loc}$ topology, 
as opposed to a pointwise Harnack estimate. 
\item{\bf ii.\ } If $\ln u$ satisfies (\ref{Eq:1:3}) for 
some $p> N+2$, then $u$ is locally analytic in the 
space variables uniformly in $t$, and $C^\infty$ in 
time. 
\end{description}
\section{Harnack Type Estimates in the Topology 
of $L^1_{\loc}$}\label{S:2}
For $\rho>0$ let $K_\rho$ be the cube centered at the origin 
of $\rn$ and edge $\rho$, and for $y\in\rn$ let $K_\rho(y)$ 
denote the homothetic cube centered at $y$. 
Moreover, $Q_\rho(\theta)$ denotes the parabolic cylinder
$K_\rho\times(-\theta\rho^2,0]$.
For $0<s<t\le T$ 
and $y\in E$ let $\rho$ be so small that
$K_{2\rho}(y)\times (s,t]\subset E_T$.
Since $u\in L^\infty_{\loc}(E_T)$ the quantity 
\begin{equation}\label{Eq:2:1}
M=\essup_{K_{2\rho}\times(s,t]} u 
\end{equation}
is well defined and finite. Also, if (\ref{Eq:1:3}) 
holds then the quantity
\begin{equation}\label{Eq:2:2}
\Lm_p=\essup_{s\le\tau\le t} 
\Big(\bint_{K_{2\rho}(y)}\Big|\ln\frac{u(x,\tau)}{M}\Big|^p 
dx\Big)^{\frac1p}
\end{equation}
is well defined and finite. 
\begin{proposition}\label{Prop:2:1} 
Let $u$ be a non-negative, local, weak solution to 
(\ref{Eq:1:1}) satisfying in addition (\ref{Eq:1:2}) and 
(\ref{Eq:1:3}) for some $p\ge2$.  There 
exists a positive constant $\gm$ depending only on 
$\{N,r,p\}$ and $\Lm_1$ and $\Lm_2$, such that 
for all cylinders $K_{2\rho}(y)\times[s,t]\subset E_T$, there holds
\begin{equation}\label{Eq:2:3}
\sup_{s<\tau<t}\int_{K_\rho(y)} u(x,\tau)dx\le
\gm\left(\inf_{s<\tau<t}\int_{K_{2\rho}(y)}u(x,\tau)dx
+\frac{t-s}{\rho^{\lm}}\right),
\end{equation}
where 
\begin{equation}\label{Eq:2:4}
\lm=2-N.
\end{equation}
\end{proposition}
\subsection{Weak Solutions Versus Distributional 
Solutions}\label{S:2:1}
The $L^1_{\loc}$ Harnack type estimate (\ref{Eq:2:3}), 
continues to hold for merely distributional solutions 
to the second of (\ref{Eq:1:1}) whereby $\ln u$ is only 
in $L^1_{\loc}(E_T)$. The assumption 
(\ref{Eq:1:2}) however is still in force, and (\ref{Eq:1:3}) 
is required to hold only for some $p>1$. The constant $\gm$ 
depends only on $\Lm_p$ for some $p>1$. 

In \S~\ref{S:Dis} we will prove (\ref{Eq:2:3}) first 
for such distributional solutions. The proof is rather 
simple due to the linearity of the principal part with 
respect to $\ln u$. The linearity however is immaterial, 
as (\ref{Eq:2:3}) is a structural inequality valid for 
weak solutions to quasilinear parabolic equations with 
singularity and degeneracy of the same nature as (\ref{Eq:1:1}).  
To be specific, consider non-negative, local, weak 
solutions to quasilinear parabolic equations of the type 
\begin{equation}\label{Eq:2:5}
\begin{aligned}
&u\in C_{\loc}\big(0,T; L^2_{\loc}(E)\big),\quad
\ln u\in L_{\loc}^2\big(0,T; W_{\loc}^{1,2}(E)\big)\\
&u_t-\dvg{\bl A}(x,t,u,Du) = 0\quad
\text{ weakly in }\> E_T.
\end{aligned}
\end{equation}
Here the function $\bl{A}:E_T\times\rr^{N+1}\to\rn$ 
is only assumed to be measurable
and subject to the structure conditions 
\begin{equation}\label{Eq:2:6}
\begin{aligned}
\bl{A}(x,t,u,p)\cdot p&\ge C_o\frac{|p|^2}u\\
|\bl{A}(x,t,u,p)|&\le C_1\frac{|p|}u
\end{aligned}\qquad \text{ a.e. in }\> E_T,
\end{equation}
where $C_o$ and $C_1$ are given positive constants. Assume 
that $u$ and $\ln u$ satisfy (\ref{Eq:1:2})--(\ref{Eq:1:3}) 
and introduce $M$ and $\Lm_p$ as in (\ref{Eq:2:1}) and 
(\ref{Eq:2:2}).  Then $u$ satisfies (\ref{Eq:2:3}) 
with $\gm$ depending on the data $\{N,r,p,C_o,C_1\}$ 
and $\Lm_1$ and $\Lm_2$. 

The proof of this fact is more involved and it is given in 
\S~\ref{S:Weak}.
\section{Analyticity of Local Weak 
Solutions to (\ref{Eq:1:1})}\label{S:3}
The precise statement of these results hinges on the 
notion of ``intrinsic neighborhood'' of a point $\pto$, as 
determined by the degeneracy and singularity of the 
equation in (\ref{Eq:1:1}). 
\subsection{The Intrinsic Geometry of (\ref{Eq:1:1}) and 
Main Results from \cite{DBGL1}}\label{S:3:1}
Let $u$ be a non-negative, local, weak solution to 
(\ref{Eq:1:1}). Having fixed $\pto\in E_T$, and 
$K_{8\rho}(x_o)\subset E$, introduce the quantity 
\begin{equation}\label{Eq:3:1} 
\theta\df{=}\varep\uqoonrm
\end{equation}
where $\varep\in(0,1)$ is to be chosen, and $q>1$ 
is arbitrary.  If $\theta>0$ assume that 
\begin{equation}\label{Eq:3:2} 
\pto+Q_{8\rho}(\theta)=
K_{8\rho}(x_o)\times(t_o-\theta(8\rho)^2, t_o]\subset E_T.
\end{equation}
These are backward, parabolic cylinders with ``vertex`` 
at $\pto$ whose height depends on the solution itself 
through the quantity $\theta$. In this sense they are 
intrinsic to the solution itself. 

Continue to assume that $u$ satisfies (\ref{Eq:1:2}) 
and let $\ln u$ satisfy (\ref{Eq:1:3}) for some $p>N+2$. 
Then 
\begin{equation}\label{Eq:3:3}
M=\essup_{\pto+Q_{8\rho}(\theta)} u 
\end{equation}
is well defined and finite. Moreover the 
dimensionless quantity
\begin{equation}\label{Eq:3:4}
\eta= \Big[\bint_{K_\rho(x_o)}\Big(\frac{u(x,t_o)}{M}\Big)^qdx
\Big]^{\frac1q\,\frac2{2r-N}}=
\Big(\frac{\theta}{\varep M}\Big)^{\frac2{2r-N}}
\end{equation}
is well defined and strictly positive. Finally 
the dimensionless quantity
\begin{equation}\label{Eq:3:5}
\Lm_p=\essup_{t_o-\theta(8\rho)^2<\tau<t_o} 
\Big(\bint_{K_{8\rho}(x_o)} \Big|\ln\frac{u(x,\tau)}{M}\Big|^p 
dx\Big)^{\frac1p},
\quad\text{ for some}\> p>N+2
\end{equation}
is well defined and finite.  
\begin{theorem}[Pointwise Harnack Estimate 
\cite{DBGL1}]\label{Thm:3:1}
Let $u$ be a non-negative, local, weak solution to (\ref{Eq:1:1}), 
satisfying the integrability conditions (\ref{Eq:1:2}) 
and (\ref{Eq:1:3}) for some $p>N+2$, and assume $\theta>0$. 
There exist a constant $\varep\in (0,1)$,  and a continuous, 
increasing function $\eta\to f(\eta,\Lm_p)$ defined in $\rr^+$ 
and vanishing at $\eta=0$, that can be quantitatively 
determined apriori only in terms of $\{N,p,q\}$, and $\Lm_p$, 
such that 
\begin{equation}\label{Eq:3:6}
\begin{aligned}
&\inf_{K_{4\rho}(x_o)} u(\cdot,t)\ge f(\eta,\Lm_p) 
\sup_{\pto+Q_{2\rho}(\frac14\theta)}u\\
{}\\
&\quad \text{for all }\>
t\in (t_o-{\txty\frac1{16}}\theta\rho^2\,,\,t_o]
\end{aligned}
\end{equation}
For $\eta\to0$ and $\Lm_p\to\infty$, the function 
$\eta\to f(\eta,\Lm_p)$ can be taken to be of the form
\begin{equation}\label{Eq:3:7}
f(\eta)= \exp\Big\{-\frac{\Lm_p^{C_1}}{\eta^{C_2}}\Big\}
\qquad\text{ for }\> 0\le \eta\ll1
\quad\text{ and }\>\Lm_p\gg1 
\end{equation}
for positive constants $C_1$ and $C_2$ that can be determined 
apriori only in terms of  $\{N,p,q\}$. 
Moreover 
\begin{equation}\label{Eq:3:8}
\varep\to0\>\text{ and }\>C_1+C_2 \to \infty
\quad\text{ as }\> p\to\infty\quad\text{ or }\> p\to N+2.
\end{equation} 
\end{theorem}
\begin{remark}\label{Rmk:3:1}
{\normalfont In \cite{DBGL1} the constant $\eta$ was given a 
more general form. For the purpose of this note the definition 
(\ref{Eq:3:4}) represents  the degeneracy of the equation, 
quantified  by $M\to\infty$. The occurrence $\Lm_p\to\infty$ 
quantifies, roughly speaking, the singularity of the equation.
}
\end{remark}
\subsection{Analyticity in the Space Variable, of 
Solutions to (\ref{Eq:1:1}) at $\pto$}\label{S:3:2}
\begin{theorem}\label{Thm:3:2}
Let $u$ be a non-negative, local, weak solution to (\ref{Eq:1:1}), 
satisfying the integrability conditions (\ref{Eq:1:2}) 
and (\ref{Eq:1:3}) for some $p>N+2$, and assume $\theta>0$. 
There exist two parameters $C$ and $H$, that have a 
polynomial dependence on $f(\eta)$, $[f(\eta)]^{-1}$, $N$, 
such that for every $N$-dimensional multi-index $\al$
\begin{equation}\label{Eq:3:9}
|D^{\al}u\pto|\le \frac{C H^{|\al|} |\al|!}{\rho^{|\al|}}\,u\pto.
\end{equation}
Moreover, for every non-negative integer $k$
\begin{equation}\label{Eq:3:10}
\Big|\frac{\pl^k}{\pl t^k}u\pto\Big|\le 
\frac{C H^{2k}(2k)!}{\rho^{2k}}\,u\pto^{1-k}.
\end{equation}
\end{theorem}
\begin{remark}\label{Rmk:4:1} {\normalfont The 
theorem continues to hold, with the same assumptions, 
for local, weak solutions to the quasilinear equations  
(\ref{Eq:2:5}), provided the function $\bl{A}$ is 
analytic in all its arguments whenever $u$ is bounded 
above and below by positive constants.
}
\end{remark}
\section{Approximating (\ref{Eq:1:1}) by Porous Medium Type 
Equations}\label{S:4}
Consider local, non-negative, weak solutions in $E_T$ to 
the porous medium equation 
\begin{equation}\label{Eq:4:1}
\begin{aligned}
& u\in C_{\loc}\big(0,T; L^2_{\loc}(E)\big),\quad 
w\in L^2_{\loc}\big(0,T;W^{1,2}_{\loc}(E)\big);\\
&u_t-\Delta w=0 \qquad\text{ weakly in }\> E_T=E\times(0,T]
\end{aligned}
\end{equation}
where 
\begin{equation}\label{Eq:4:2} 
w=\frac{u^{m}-1}m\qquad\text{ for }\> 0<m\ll1.
\end{equation}
As $m\to0$, formally (\ref{Eq:4:1})--(\ref{Eq:4:2}) tend to 
(\ref{Eq:1:1}). In \cite{DBGL2} a precise topology was 
introduced by which such a formal limit is rigorous. A natural 
question is whether solutions to (\ref{Eq:4:1})--(\ref{Eq:4:2}) 
satisfy a version of the $L^1_{\loc}$ Harnack estimate 
(\ref{Eq:2:3}), which as $m\to0$ tends, in some appropriate 
sense to be made precise, to that of Proposition~\ref{Prop:2:1}. 
A similar issue arises for the local analyticity of 
Theorem~\ref{Thm:3:2}. 
\subsection{Harnack Type Estimates in the Topology 
of $L^1_{\loc}$, for Weak Solutions to 
(\ref{Eq:4:1})--(\ref{Eq:4:2})}\label{S:4:1}

A first statement in this direction is that $u$ satisfies 
\begin{equation}\label{Eq:4:3}
\sup_{s<\tau<t}\int_{K_\rho(y)} u(x,\tau)dx\le
\gm\left[\inf_{s<\tau<t}\int_{K_{2\rho}(y)}u(x,\tau)dx
+\Big(\frac{t-s}{\rho^{\lm}}\Big)^{\frac1{1-m}}\right],
\end{equation}
where 
\begin{equation}\label{Eq:4:4}
\lm=N(m-1)+2.
\end{equation}
Here $\gm$ depends upon $N$ and $m$ and $\gm(m)\to\infty$ 
as $m\to0$. Thus, one cannot formally recover (\ref{Eq:2:3}) 
by letting $m\to0$ in (\ref{Eq:4:3}). However, (\ref{Eq:4:3}) 
is rather general as it continues to hold for non-negative, 
local weak solutions to general quasi-linear version 
of (\ref{Eq:4:1}). Precisely 
\begin{equation}\label{Eq:4:5}
\begin{aligned}
&u\in C_{\loc}\big(0,T; L^2_{\loc}(E)\big),\quad
w\in L_{\loc}^2\big(0,T; W_{\loc}^{1,2}(E)\big)\\
&u_t-\dvg{\bl A}(x,t,u,Du) = 0\quad
\text{ weakly in }\> E_T
\end{aligned}
\end{equation}
where the function $\bl{A}:E_T\times\rr^{N+1}\to\rn$ 
is only assumed to be measurable
and subject to the structure conditions 
\begin{equation}\label{Eq:4:6}
\begin{aligned}
\bl{A}(x,t,u,p)\cdot p&\ge C_o{u^{m-1}}{|p|^2}\\
|\bl{A}(x,t,u,p)|&\le C_1{u^{m-1}}{|p|}
\end{aligned} \qquad\text{ a.e. in }\> E_T,
\end{equation}
where $C_o$ and $C_1$ are given positive constants. In such 
a case the constant $\gm$ in (\ref{Eq:4:3}) depends also 
on these structural constants. The proof of these statements 
is in \cite{DBGV-mono}, Appendix~B.

A major difference between (\ref{Eq:2:3}) and (\ref{Eq:4:3}) 
is that in the latter $u$ is not required to be locally 
bounded, nor does $\gm$ depend on some analogue of the quantity 
$\Lm_p$ as defined in (\ref{Eq:2:2}). This raises the question 
as to whether (\ref{Eq:4:3}) holds with $\gm$ independent 
of $m$ but dependent on some analogue of $\Lm_p$.
 
Henceforth we assume 
\begin{equation}\label{Eq:4:7}
u\in L^r_{\loc}(E_T)\quad\text{ for some } 
\> r>\max\big\{1;{\txty\frac{N}{2}(1-m)}\big\}
\end{equation}
and that 
\begin{equation}\label{Eq:4:8}
w\in L^\infty_{\loc}\big(0,T;L^p_{\loc}(E)\big)
\quad\text{ for some } \> p\ge 1.
\end{equation}  
It was shown in \cite{DBGV-mono} that (\ref{Eq:4:7}) 
implies that $u\in L^\infty_{\loc}(E_T)$ 
and hence the corresponding quantity $M$ defined 
as in (\ref{Eq:2:1}) is well defined and finite. Set 
\begin{equation}\label{Eq:4:9}
\Lm_{m,p} =\essup_{s<\tau<t} \left(\int_{K_{2\rho}(y)} 
\Big(\frac{M^m-u(x,\tau)^m}{m\,M^m}\Big)^p dx\right)^{\frac1p}.
\end{equation}
This is the analogue of (\ref{Eq:2:2}), and, if 
(\ref{Eq:4:8}) holds, it  is well defined and finite.  
\begin{proposition}\label{Prop:4:1} Let $u$ be a 
non-negative, local, weak solution to the singular 
equations (\ref{Eq:4:5})--(\ref{Eq:4:6}), in $E_T$.  There 
exists a positive constant $\gm$ depending only on 
the data $N,C_o,C_1$,  and $\Lm_{\frac m2,1}$, 
$\Lm_{\frac m2,2}$ and independent of $m$, such 
that (\ref{Eq:4:3})--(\ref{Eq:4:4}) holds true, 
for all cylinders $K_{2\rho}(y)\times[s,t]\subset E_T$.
\end{proposition}
As a consequence, (\ref{Eq:2:3}) can be recovered from 
(\ref{Eq:4:3}), with the indicated dependences, as $m\to0$, 
provided proper conditions are placed, that insure the 
pointwise convergence of the solutions to
(\ref{Eq:4:1})-(\ref{Eq:4:2}) to solutions of (\ref{Eq:1:1}). 
These conditions are identified in \cite{DBGL2} and we will 
touch on them briefly in the next subsection.
\subsection{Analyticity in the Space Variable, of 
Solutions to (\ref{Eq:4:1}) at $\pto$}\label{S:4:2}
Having fixed $\pto\in E_T$ and $K_{8\rho}(x_o)\subset E$, 
the intrinsic geometry of (\ref{Eq:4:1})--(\ref{Eq:4:2}) 
is determined by
\begin{equation}\label{Eq:4:10}
\theta_m=\varep\Big(\bint_{K_{\rho}(x_o)} 
u^q(x,t_o) dx\Big)^{\frac{1-m}q}.
\end{equation}
The intrinsic cylinders are as in (\ref{Eq:3:2}) with $\theta$ 
replaced by $\theta_m$. The analogues of $\eta$ in (\ref{Eq:3:4}) 
are 
\begin{equation}\label{Eq:4:11}
\sig=\Big[\bint_{K_\rho(x_o)}\big(\frac{u(x,t_o)}{M}\Big)^qdx 
\Big]^{\frac1q\,\frac{2}{\lm_r}}
\end{equation}
where $r\ge1$ is any number such that
\begin{equation}\label{Eq:4:12}
\lm_r=N(m-1)+2r>0.
\end{equation}
In \cite{DBGV-mono} a Harnack estimate of the form 
of (\ref{Eq:3:6}) was proved for these solutions with 
$f(\cdot)$ depending only on $\sig$ and of the form 
\begin{equation}\label{Eq:4:13}
f(\sig)=\frac{\sig^\be}{\gm(m)}
\end{equation} 
where $\gm(m)\to\infty$ as $m\to0$. 
The constant $\be$ depends on $\lm_r$ and $\be(\lm_r)\to\infty$ 
as $\lm_r\to0$. It was observed in \cite{DBGV-mono} \S~21.5.3 
that, for each fixed $m\in(0,1)$,  such an estimate implies 
the local analyticity of the solutions in the space 
variables about $\pto$, and at least the Lipschitz continuity 
in time. 

Because of the indicated dependence of $\gm(m)$ on $m$ in 
(\ref{Eq:4:13}), such a regularity does not directly carry to 
the limit as $m\to0$. In \cite{DBGL2} we established 
a Harnack estimate of the form (\ref{Eq:3:6}) for solutions 
to (\ref{Eq:4:1})--(\ref{Eq:4:2}) and its quasi-linear versions 
(\ref{Eq:4:5})-(\ref{Eq:4:6}), with $f(\cdot)$ 
depending on $\sig$, as defined in (\ref{Eq:4:11}), and 
$\Lm_{m,p}$ as defined in (\ref{Eq:4:9}) provided $p>N+2$. 
The form of such $f(\cdot)$ is the same as that 
in (\ref{Eq:3:7}) with the proper change in symbolism. 
The new feature of such an $f(\cdot)$ is that, while 
depending on the quantities $\sig$ and $\Lm_{m,p}$, each 
quantifying the degeneracy and the singularity of 
the equation, is independent of $m$ and hence is ``stable'' 
as $m\to0$, provided $\sig$ and $\Lm_{m,p}$ are 
uniformly estimated with respect to $m$.  

As a consequence, the analyticity estimates 
of Theorem~\ref{Thm:3:2}, can be recovered from the 
analogous ones for solutions to (\ref{Eq:4:1})--(\ref{Eq:4:2}) 
whenever solutions $\{u_m\}$ to the latter converge pointwise 
to solutions to the former. In \cite{DBGL2} it was shown 
that this occurs if there exists $m_{**}\in(0,1)$ such that 
\begin{equation}\label{Eq:4:14}
\begin{array}{ll} 
u_m\in L^\infty_{\loc}\big(0,T;L^r_{\loc}(E)\big) \quad
&\text{for some $r>\max\{1;\frac12 N\}$}\\
w_m\in L^\infty_{\loc}(0,T; L^p_{\loc}(E)\big)&\text{for some $p>N+2$}
\end{array}
\end{equation}
uniformly in $m\in(0,m_{**})$. It is also required that 
there exists an open set $E_o\subset E$ and a positive number 
$\sig_{E_o;T}$ such that 
\begin{equation}\label{Eq:1:15}
\int_{E_o} u_m(\cdot,T)dx\ge \sig_{E_o;T}\quad
\text{ uniformly in }\> m.
\end{equation}
\section{Proof of Proposition~\ref{Prop:2:1} for Distributional 
Solutions to (\ref{Eq:1:1})}\label{S:Dis}
The proof is a local version of an argument of 
\cite{HP} for global solutions to the porous 
medium equation for $0<m<1$.  Let $\z_1\in C^\infty_o(\rn)$ 
be such that
\begin{equation}\label{Eq:Dis:1}
\left\{
\begin{array}{l}
0\le\z_1\le1\\
\z_1=1\quad\text{ in }\quad K_\rho\\
\z_1=0\quad\text{ in }\quad \rn\backslash K_{2\rho}.
\end{array}\right .
\end{equation}
Then, by the divergence theorem
\begin{equation}\label{Eq:Dis:2}
\int_{\rn}\Dl\z_1 dx=\int_{K_{2\rho}}\Dl\z_1\,dx
=\int_{\partial K_{2\rho}}\frac{\partial\z_1}{\partial n}\,ds=0.
\end{equation}
By (\ref{Eq:Dis:2}), for any positive constant $M$, 
any $\z_1$ as in (\ref{Eq:Dis:1}), and any 
non-negative function $v$ such that $\Dl\z_1\ln v$ 
is integrable, we have
\begin{equation}\label{Eq:Dis:3}
\int_{K_{2\rho}}\Dl\z_1\ln v dx
=\int_{K_{2\rho}}\Dl\z_1\ln (\frac vM) dx.
\end{equation}
Now consider $\z_2\in C^\infty_o(0,+\infty)$ and 
$\z_1$ as in (\ref{Eq:Dis:1}), such that
\begin{equation*}
|D\z_1|\le\frac{C_1(N)}\rho,\qquad|\Dl\z_1|
\le\frac{C_2(N)}{\rho^2}.
\end{equation*}
By the previous notation, with $u$ a solution 
to (\ref{Eq:1:1}), we have
\begin{equation*}
-\int_0^\infty\int_{\rn}\z_2^{\prime}\z_1 u dxdt
=\int_0^\infty\int_{\rn}\z_2\Dl\z_1\ln u dxdt.
\end{equation*}
Taking into account (\ref{Eq:Dis:3}), for any positive 
constant $M$
\begin{equation*}
\frac d{dt}\int_{K_{2\rho}}\z_1 u dx
=\int_{K_{2\rho}}\Dl\z_1\ln(\frac uM) dx\qquad\text{ in }
\qquad{\cal D}^{\prime}(0,T),
\end{equation*}
and also in $L^1_{loc}(0,T)$. Therefore
\begin{equation*}
\Big|\frac d{dt}\int_{K_{2\rho}}\z_1 u dx\Big|
=\Big|\int_{K_{2\rho}}\Dl\z_1\ln(\frac uM)dx\Big|
\le\int_{K_{2\rho}}|\Dl\z_1| \Big|\ln(\frac uM)\Big|dx.
\end{equation*}
By the definition of $\Lm_p$ and from the previous 
estimate, 
\begin{align*}
\Big|\frac d{dt}\int_{K_{2\rho}}\z_1 u dx\Big|&\le\Lm_p
\left(\int_{K_{2\rho}}|\Dl\z_1|^{p^{\prime}}dx
\right)^{\frac1{p^{\prime}}}|K_{8\rho}|^{\frac1p}\\ 
&\le\Lm_p\frac{C(N)}{\rho^2}\rho^{\frac N{p^{\prime}}}\rho^{\frac Np}
=\frac{C(\Lm_p,N)}{\rho^\lm},
\end{align*}
where $\lm=2-N$. Taking into account
the size of the support of $\z_1$, for any $0<s<t<T$ we conclude
\begin{equation*}
\int_{K_\rho}u(x,t)dx\le C(\Lm_p,N)
\left(\int_{K_{2\rho}}u(x,s) dx+\frac{t-s}{\rho^{\lm}}\right).
\tag*{\bbox}
\end{equation*}
\section{Proof of Proposition~\ref{Prop:2:1} for Weak 
Solutions to (\ref{Eq:2:5})--(\ref{Eq:2:6})}\label{S:Weak}
\subsection{An Auxiliary Lemma}\label{S:Weak:1}
\begin{lemma}\label{Lm:Weak:1} 
Let $u$ be a non-negative, local, weak solution to 
the quasi-linear singular equations 
(\ref{Eq:2:5})--(\ref{Eq:2:6}), in $E_T$.  There exist
two positive constants $\gm_1$, $\gm_2$ depending only on the 
data $\{N,C_o,C_1\}$, such that for all cylinders 
$K_{4\rho}(y)\times[s,t]\subset E_T$, and all 
$\sig\in(0,1)$, 
\begin{equation*}
\int_s^t\int_{K_\rho(y)}\frac{|Du|^2}{u^2}\z^2dx d\tau
\le\gm_1(1+\Lm_1)\mcl{S}_\sig
+\frac{\gm_2}{\sig^2}(\Lm_1^2+\Lm_2^2)
\left(\frac{t-s}{\rho^\lm}\right),
\end{equation*}
where 
\begin{equation*}
\mcl{S}_\sig=\sup_{s <\tau<t}\int_{K_{(1+\sig)\rho}(y)} 
u(\cdot,\tau)dx.
\end{equation*}
\end{lemma}
\noi{\it Proof -} 
Assume $(y,s)=(0,0)$, fix $\sig\in(0,1)$,
and let $x\to\z(x)$ be a non-negative piecewise smooth 
cutoff function in $K_{(1+\sig)\rho}$ that vanishes outside 
$K_{(1+\sig)\rho}$, equals one on $K_\rho$, and such that 
$|D\z|\le(\sig\rho)^{-1}$.  Let $s_1\in[0,t]$ be such that
\begin{equation*}
\ssig =\sup_{0 <s<t}\int_{K_{(1+\sig)\rho}(y)} 
u(\cdot,s)dx=\int_{K_{(1+\sig)\rho}(y)} 
u(\cdot,s_1)dx.
\end{equation*}
We also set
\begin{equation*}
\ssigbar \df=\frac{\mcl{S}_\sig}{\rho^N}.
\end{equation*}
In the weak formulation of (\ref{Eq:2:5})--(\ref{Eq:2:6}), 
take the test function 
\begin{equation*}
\vp=\left(\ln\frac{\ssigbar}u\right)_+ \z^2
\end{equation*}
and integrate over $Q=K_{(1+\sig)\rho}\times(0,t]$, 
to obtain
\begin{align*}
0=&\iint_{Q}\frac{\pl}{\pl\tau} u 
\left(\ln\frac{\ssigbar}u\right)_+\z^2dx\,d\tau
+\iint_{Q} \bl{A}(x,\tau,u,Du)\cdot 
D\left[\left(\ln\frac{\ssigbar}u\right)_+\z^2\right] dx\,d\tau\\
&= I_1+I_2.
\end{align*}
We estimate these two terms separately. 
\begin{equation*}
\begin{aligned}
I_1&=\iint_{Q} \frac{\pl}{\pl\tau} u 
\left(\ln\frac{\ssigbar}u\right)_+\z^2dx\,d\tau
=\iint_{Q\cap[u<\ssigbar]} \frac{\pl}{\pl\tau}u 
\left(\ln\frac{\ssigbar}u\right)\z^2dx\,d\tau\\
&=\int_{K_{(1+\sig)\rho}\cap[u<\ssigbar]}\z^2(x)
\left(u\ln\frac{\ssigbar}u+u\right)(x,t)\,dx\\
&-\int_{K_{(1+\sig)\rho}\cap[u<\ssigbar]}\z^2(x)
\left(u\ln\frac{\ssigbar}u+u\right)(x,0)\,dx.
\end{aligned}
\end{equation*}
Next,
\begin{align*} 
I_2&=\iint_{Q} \bl{A}(x,\tau,u,Du)\cdot 
D\left[\left(\ln\frac{\ssigbar}u\right)_+\z^2\right] dx\,d\tau\\
&=\iint_{Q\cap[u<\ssigbar]} \bl{A}(x,\tau,u,Du)\cdot 
D\left[\left(\ln\frac{\ssigbar}u\right)\z^2\right] dx\,d\tau\\
&=-\iint_{Q\cap[u<\ssigbar]}\z^2
\bl{A}(x,\tau,u,Du)\frac{Du}u dx\, d\tau\\
&+2\iint_{Q\cap[u<\ssigbar]}\z\,\ln\frac{\ssigbar}u
\bl{A}(x,\tau,u,Du){D\z} dx\, d\tau\\
&\le-C_o\iint_{Q\cap[u<\ssigbar]}\z^2
\frac{|Du|^2}{u^2}dx\,d\tau\\
&+2C_1\iint_{Q\cap[u<\ssigbar]}\z\,
\ln\frac{\ssigbar}u\frac{|Du|}u|D\z| dx\, d\tau\\
&\le-\frac{C_o}2\iint_{Q\cap[u<\ssigbar]}
\z^2\frac{|Du|^2}{u^2} dx\, d\tau +
\frac{\gm}{\sig^2\rho^2}\iint_{Q\cap[u<\ssigbar]}
\left|\ln\frac{\ssigbar}u\right|^2 dx\,d\tau,
\end{align*}
where $\gm=2\frac{C_1^2}{C_o}$.
Therefore, we conclude that
\begin{align*}
&\frac{C_o}2\iint_{Q\cap[u<\ssigbar]}\z^2
\frac{|Du|^2}{u^2} dx\, d\tau\le
\int_{K_{(1+\sig)\rho}\cap[u<\ssigbar]}\z^2(x)
\left(u\ln\frac{\ssigbar}u+u\right)(x,t)\,dx\\
&+\frac{\gm}{\sig^2\rho^2}\iint_{Q\cap[u<\ssigbar]}
\left|\ln\frac{\ssigbar}u\right|^2 dx\,d\tau\\
&=J_1+J_2.
\end{align*}
We have
\begin{align*}
J_1&=\int_{K_{(1+\sig)\rho}\cap[u<\ssigbar]}
\z^2(x)\left(u\ln\frac{\ssigbar}u+u\right)(x,t)\,dx\\
&\le\int_{K_{(1+\sig)\rho}\cap[u<\ssigbar]}
\left(u\ln\frac{\ssigbar}u+u\right)(x,t)\,dx\\
&=\ssigbar\int_{K_{(1+\sig)\rho}\cap[u<\ssigbar]}
\frac u{\ssigbar}\ln\frac{\ssigbar}u(x,t)\,dx
+\int_{K_{(1+\sig)\rho}\cap[u<\ssigbar]}u(x,t)\,dx\\
&\le\ssigbar\int_{K_{(1+\sig)\rho}\cap[u<\ssigbar]}
\ln\frac{\ssigbar}u(x,t)\,dx+\ssig\\
&\le\gm\ssig\bint_{K_{(1+\sig)\rho}}
\ln \frac Mu(x,t)\,dx+\ssig
\le\gm\Lm_1\ssig+\ssig=\gm(1+\Lm_1)\ssig,
\end{align*}
where $\gm=2^N$. Moreover,
\begin{align*}
J_2&=\frac{\gm}{\sig^2\rho^2}\iint_{Q\cap[u<\ssigbar]}
\left|\ln\frac{\ssigbar}u\right|^2 dx\,d\tau
\le \frac{\gm}{\sig^2\rho^2}\iint_{Q\cap[u<\ssigbar]}
\left|\ln\frac Mu\right|^2 dx\,d\tau\\
&\le \frac{\gm}{\sig^2\rho^2}\iint_{Q}
\left|\ln\frac Mu\right|^2 dx\,d\tau=\frac{\gm}{\sig^2\rho^2}
\int_0^t\int_{K_{(1+\sig)\rho}}
\left|\ln\frac Mu\right|^2 dx\,d\tau\\
&\le\frac{\gm}{\sig^2\rho^2}\sup_{0<\tau<t}
\bint_{K_{(1+\sig)\rho}}\left|\ln\frac Mu\right|^2 dx\,t\rho^N
\le\frac{\gm}{\sig^2}\Lm_2^2\left(\frac t{\rho^\lm}\right).
\end{align*}
Hence, we have
\begin{equation}\label{Eq:Weak:1}
\iint_{Q\cap[u<\ssigbar]}\z^2\frac{|Du|^2}{u^2} dx\, d\tau
\le\gm(\Lm_1+1)\ssig+\frac{\gm}{\sig^2}\Lm_2^2
\left(\frac t{\rho^\lm}\right).
\end{equation}
Now, if we take the test function 
\begin{equation*}
\vp=\left(\ln\frac u{\ssigbar}\right)_+ \z^2
\end{equation*}
in the weak formulation of (\ref{Eq:2:5})--(\ref{Eq:2:6}) 
and integrate over $Q=K_{(1+\sig)\rho}\times(0,t]$, 
we obtain
\begin{align*}
0=&\iint_{Q}\frac{\pl}{\pl\tau} u 
\left(\ln\frac u{\ssigbar}\right)_+\z^2dx\,d\tau
+\iint_{Q} \bl{A}(x,\tau,u,Du)\cdot 
D\left[\left(\ln\frac u{\ssigbar}\right)_+\z^2\right] dx\,d\tau\\
&= I_3+I_4.
\end{align*}
We estimate these two terms separately. 
\begin{equation*}
\begin{aligned}
I_3&=\iint_{Q} \frac{\pl}{\pl\tau} u 
\left(\ln\frac u{\ssigbar}\right)_+\z^2dx\,d\tau
=\iint_{Q\cap[u>\ssigbar]} \frac{\pl}{\pl\tau} u 
\left(\ln\frac u{\ssigbar}\right)\z^2dx\,d\tau\\
&=\int_{K_{(1+\sig)\rho}\cap[u>\ssigbar]}\z^2(x)
\left(u\ln\frac u{\ssigbar}-u\right)(x,t)\,dx\\
&-\int_{K_{(1+\sig)\rho}\cap[u>\ssigbar]}\z^2(x)
\left(u\ln\frac u{\ssigbar}-u\right)(x,0)\,dx.
\end{aligned}
\end{equation*}
Next,
\begin{align*} 
I_4&=\iint_{Q} \bl{A}(x,\tau,u,Du)\cdot 
D\left[\left(\ln\frac u{\ssigbar}\right)_+\z^2\right] dx\,d\tau\\
&=\iint_{Q\cap[u>\ssigbar]} \bl{A}(x,\tau,u,Du)\cdot 
D\left[\left(\ln\frac u{\ssigbar}\right)\z^2\right] dx\,d\tau\\
&=\iint_{Q\cap[u>\ssigbar]}\z^2\bl{A}(x,\tau,u,Du)
\frac{Du}u dx\, d\tau\\
&+2\iint_{Q\cap[u>\ssigbar]}\z\,\ln\frac u{\ssigbar}
\bl{A}(x,\tau,u,Du){D\z} dx\, d\tau\\
&\ge C_o\iint_{Q\cap[u>\ssigbar]}\z^2
\frac{|Du|^2}{u^2}dx\,d\tau\\
&-2C_1\iint_{Q\cap[u>\ssigbar]}\z\,
\ln\frac u{\ssigbar}\frac{|Du|}u|D\z| dx\, d\tau\\
&\ge\frac{C_o}2\iint_{Q\cap[u>\ssigbar]}\z^2
\frac{|Du|^2}{u^2} dx\, d\tau -
\frac{\gm}{\sig^2\rho^2}\iint_{Q\cap[u>\ssigbar]}
\left|\ln\frac u{\ssigbar}\right|^2 dx\,d\tau,
\end{align*}
where again $\gm=2\frac{C_1^2}{C_o}$. Therefore, we conclude that
\begin{align*}
&\frac{C_o}2\iint_{Q\cap[u>\ssigbar]}\z^2
\frac{|Du|^2}{u^2} dx\, d\tau
\le\int_{K_{(1+\sig)\rho}\cap[u>\ssigbar]}\z^2(x)
\left(u\ln\frac u{\ssigbar}\right)(x,0)\,dx\\
&+\int_{K_{(1+\sig)\rho}\cap[u>\ssigbar]}\z^2(x) u(x,t)\,dx+
\frac{\gm}{\sig^2\rho^2}\iint_{Q\cap[u>\ssigbar]}
\left|\ln\frac u{\ssigbar}\right|^2 dx\,d\tau\\
&\le\int_{K_{(1+\sig)\rho}\cap[u>\ssigbar]}\z^2(x)
\left(u\ln\frac u{\ssigbar}\right)(x,0)\,dx
+\ssig\\
&+\frac{\gm}{\sig^2\rho^2}\iint_{Q\cap[u>\ssigbar]}
\left|\ln\frac u{\ssigbar}\right|^2 dx\,d\tau\\
&\le\int_{K_{(1+\sig)\rho}\cap[u>\ssigbar]}\z^2(x)
\left(u\ln\frac M{\ssigbar}\right)(x,0)\,dx
+\ssig+\frac{\gm}{\sig^2\rho^2}
\left|\ln\frac M{\ssigbar}\right|^2\,t\rho^N\\
&\le\ssig+\left(\ln\frac M{\ssigbar}\right)\ssig
+\frac{\gm}{\sig^2}\left|\ln\frac M{\ssigbar}\right|^2
\left(\frac t{\rho^\lm}\right).
\end{align*}
We need to evaluate $\dsty\ln M/{\ssigbar}$. As in 
the interval $(0,1]$ the function $f(s)=-\ln s$ 
is convex, Jensen's inequality yields
\begin{align*}
&\ln\frac M{\ssigbar}=\ln\frac M{\bint_{K_{(1+\sig)\rho}}u(x,s_1)dx}
=-\ln \bint_{K_{(1+\sig)\rho}}\frac{u(x,s_1)}Mdx\\
&\le\bint_{K_{(1+\sig)\rho}}-\ln\frac{u(x,s_1)}Mdx
=\bint_{K_{(1+\sig)\rho}}\ln\frac M{u(x,s_1)}dx\\
&\le\gm\Lm_1,
\end{align*}
where $\gm=2^N$. Hence, we have
\begin{equation}\label{Eq:Weak:2}
\iint_{Q\cap[u>\ssigbar]}\z^2\frac{|Du|^2}{u^2} dx\, d\tau
\le\gm(\Lm_1+1)\ssig+\frac{\gm}{\sig^2}\Lm_1^2
\left(\frac t{\rho^\lm}\right).
\end{equation}
The lemma follows by combining estimates 
(\ref{Eq:Weak:1}) and (\ref{Eq:Weak:2}). 

The use of 
$\left(\ln\frac {\ssigbar}u\right)_+ \z^2$
as test function can be justified using  $\left(\ln\frac {\ssigbar}{u+\eps}\right)_+$ and then letting 
$\eps\to0$.\hfill\bbox
\begin{corollary}\label{Cor:Weak:1} 
Let $u$ be a non-negative, local, weak solution to 
the singular equations (\ref{Eq:2:5})--(\ref{Eq:2:6}), in $E_T$.  
There exists a positive constant $\gm$ depending only on the 
data $\{N,C_o,C_1\}$, such that for all cylinders 
$K_{4\rho}(y)\times[s,t]\subset E_T$, and all 
$\sig\in(0,1)$,  
\begin{align*}
&\frac1{\rho}
\int_s^t\int_{K_\rho(y)}|\bl{A}(x,\tau,u,Du)|dx\,d\tau\\
&\le \gm\max\left\{(1+\Lm_1)^{\frac12};(\Lm_1^2+\Lm_2^2)^{\frac12}
\right\}\left[\ssig+\frac1{\sig^2}
\Big(\frac{t-s}{\rho^\lm}\Big)\right]^{\frac12}
\left(\frac {t-s}{\rho^\lm}\right)^{\frac12}.
\end{align*}
\end{corollary}
\noi{\it Proof -} Assume $(y,s)=(0,0)$, and 
let $Q=K_\rho\times(0,t]$. By the structure conditions 
of $\bl{A}$ 
\begin{align*}
\frac1{\rho}
\int_0^t\int_{K_\rho}
|\bl{A}(x,\tau,u,Du)|dx\,d\tau
&\le \frac{C_1}{\rho}\iint_Q \frac{|Du|}udx\,d\tau\\
&\le \frac{C_1}{\rho}\left(\iint_Q 
\frac{|Du|^2}{u^2}dx\,d\tau\right)^{\frac12}
\rho^{\frac N2}t^{\frac12}\\
&=C_1\left(\iint_Q \frac{|Du|^2}{u^2}dx\,d\tau
\right)^{\frac12}\left(\frac t{\rho^\lm}\right)^{\frac12}.
\end{align*}
By Lemma~\ref{Lm:Weak:1} we conclude.\hfill\bbox
\subsection{Proof of Proposition~\ref{Prop:2:1}}\label{S:Weak:2}
Assume $(y,s)=(0,0)$. For $n=0,1,2\dots$ set
\begin{equation*}
\rho_n=\bsu_{j=1}^n\frac1{2^j}\rho,\quad K_n=K_{\rho_n};\qquad 
\trn=\frac{\rho_n+\rho_{n+1}}2,\quad 
\tkn=K_{\trn}
\end{equation*}
and let $x\to\z_n(x)$ be a non-negative, piecewise 
smooth cutoff function in $\tkn$ that equals
one on $K_n$, and such that $|D \z_n|\le 2^{n+2}/\rho$. 
In the weak formulation of (\ref{Eq:2:5})--(\ref{Eq:2:6}) 
take $\z_n$ as a test function, to obtain
\begin{align*}
\int_{\tkn}u(x,\tau_1)\z_n dx&\le\int_{\tkn} u(x,\tau_2)\z_ndx
+\frac{2^{n+2}}{\rho}
\Big|\int_{\tau_1}^{\tau_2}\int_{\tkn}
|\bl{A}(x,\tau,u,Du)|dx\,d\tau\Big|\\
&\le \int_{\tkn} u(x,\tau_2)\z_ndx+\gm\, 2^n 
\mcl{S}_{n+1}^{\frac{1}2}\Big(\frac{t}{\rho^\lm}
\Big)^{\frac12}+\gm\, 4^n\Big(\frac{t}{\rho^\lm}\Big),
\end{align*}
where 
\begin{equation*}
\mcl{S}_n= \sup_{0\le\tau\le t} \int_{K_n} u(\cdot,\tau)dx. 
\end{equation*}
Since the time levels $\tau_1$ and $\tau_2$ are 
arbitrary, choose $\tau_2$ one for which 
\begin{equation*}
\int_{K_{2\rho}} u(\cdot,\tau_2)dx =\inf_{0\le\tau\le t} 
\int_{K_{2\rho}} u(\cdot,\tau) dx\,\df{=}\, \mcl{I}. 
\end{equation*}
With this notation, the previous inequality takes the form 
\begin{equation*}
\mcl{S}_n\le \mcl{I} + 
\gm\big(\text{data},\Lm_1,\Lm_2\big)\, 4^n 
\Big(\frac{t}{\rho^\lm}\Big)
+\gm\big(\text{data},\Lm_1,\Lm_2\big)\, 2^n
\mcl{S}_{n+1}^{\frac{1}2}
\Big(\frac{t}{\rho^\lm}\Big)^{\frac12}.
\end{equation*}
By Young's inequality, for all $\varep_o\in(0,1)$ 
\begin{equation*}
\mcl{S}_n\le\varep_o\mcl{S}_{n+1} 
+\gm\big(\text{data},\Lm_1,\Lm_2,\varep_o\big) 4^n
\Big[\mcl{I}+\Big(\frac{t}{\rho^\lm}\Big)\Big].
\end{equation*}
From this, by iteration 
\begin{equation*}
\mcl{S}_o\le \varep_o^n\mcl{S}_n+
\gm(\text{data},\Lm_1,\Lm_2,\varep_o)\Big[\mcl{I}+
\Big(\frac{t}{\rho^\lm}\Big)\Big]
\bsu_{i=0}^{n-1} (4\varep_o)^i.
\end{equation*}
Choose $\varep_o$ so that the last term is majorized by a 
convergent series, and let $n\to\infty$.\hfill\bbox

The proof of Proposition~\ref{Prop:4:1} for weak 
solutions to the porous medium type quasilinear 
equations (\ref{Eq:4:5})--(\ref{Eq:4:6}), is 
similar, with the obvious modifications, and we 
confine it to Appendix~B.
\section{Analyticity in the Space Variables, of Solutions 
to (\ref{Eq:1:1})}\label{S:An:1}
Let $u$ be a non-negative, local, weak solution 
to (\ref{Eq:1:1}), satisfying the integrability 
conditions (\ref{Eq:1:2}) and (\ref{Eq:1:3}) for some $p>N+2$. 

Fix $\pto\in E_T$, assume that $K_{8\rho}(x_o)\subset E$, 
and assume that the quantity $\theta$ defined in 
(\ref{Eq:3:1}) is positive. The cylinder $\pto+Q_{8\rho}(\theta)$ 
is assumed to be contained in the domain of definition of $u$ 
as in (\ref{Eq:3:2}). The quantities $M$, $\eta$ and $\Lm_p$ 
are defined as in (\ref{Eq:3:3})--(\ref{Eq:3:5}). 

From the Harnack-type inequality \eqref{Eq:3:6}, 
\begin{equation}\label{Eq:An:1}
[f(\eta)]\,u\pto\le u(x,t)\le [f(\eta)]^{-1}u\pto
\end{equation}
for any $(x,t)$ within the cylinder
\begin{equation}\label{Eq:An:2}
Q\equiv K_{2\rho}(x_o)\times(t_o-\frac1{16}\theta\rho^2,t_o].
\end{equation}
Equation (\ref{Eq:1:1}) can be rewritten as
\begin{equation}\label{Eq:An:3}
u_t-\dvg(\frac1u\,Du)=0.
\end{equation}
By (\ref{Eq:An:1}) this can be regarded as a particular 
instance of a linear parabolic equation with bounded 
and measurable coefficients. By known results 
(for example, \cite{LSU}, Chapter~II) local, 
weak solutions to \eqref{Eq:An:3} are locally bounded 
and locally H\"older continuous. Consequently, (\ref{Eq:An:3})
can be regarded as a linear parabolic equation with 
bounded, and H\"older continuous coefficients. Again 
by classical theory (see \cite{LSU}, Chapter~III), 
one can conclude that local, weak solutions are 
indeed $C^\infty$ with respect to the space variable.

By (\ref{Eq:An:1}) the quantity $\theta$ 
can be estimated as
\begin{equation}\label{Eq:An:4}
\begin{aligned} 
&\theta\le\varep\sup_{K_\rho(x_o)}u(\cdot,t_o)
\le\varep[f(\eta)]^{-1}u\pto\\ 
&\theta\ge\varep\inf_{K_\rho(x_o)}u(\cdot,t_o)
\ge\varep[f(\eta)]\,u\pto.
\end{aligned}
\end{equation}
Let $\dl=\dl(\eta)\df=\varep[f(\eta)]$ and introduce 
the change of variables
\begin{equation*}
x\to\frac{x-x_o}\rho,\quad t\to\frac{t-t_o}{u\pto\rho^2},
\quad v=\frac{u}{u\pto}.
\end{equation*}
It maps $Q$ onto to
\begin{equation}\label{Eq:An:5}
\tilde Q\df=K_2\times(-\frac1{16}\frac{\theta}{u\pto},0]
\supset Q_\dl\df=K_2\times(-\frac\dl{16},0],
\end{equation}
and within $Q_\dl$ the function $v$ satisfies
\begin{equation}\label{Eq:An:6}
v_t-\frac1v\Delta v=-\frac{|Dv|^2}{v^2},
\end{equation}
with 
\begin{equation*}
f(\eta)\le v\le f(\eta)^{-1}.
\end{equation*} 
By a result of \cite{KN}, there exist constants $0<\sig<1$,
$C$ and $H$ such that  
\begin{equation}\label{Eq:An:7}
\begin{aligned}
&\sup_{Q_{\sig\dl}} |D^{\al}v|\le C H^{|\al|} |\al|!,\\
&\sup_{Q_{\sig\dl}} \left|\frac{\pl^k}{\pl t^k}v\right|\le C H^{2k}(2k)!
\end{aligned}
\end{equation}
where $Q_{\sig\dl}=K_{2\sig}\times(-\frac1{16}\sig\dl,0]$.
Tracing the dependence of constants gives
\begin{equation}
C=\gm_1 C_o, \quad H=\gm_2 \max\{C_o[f(\eta)]^{-1}, 
[f(\eta)]^{-2}\} 
\end{equation}
where $\gm_1$ and $\gm_2$ are constants independent 
of $v$ and $C_o$ is a function of $f(\eta)$ and satisfies
\begin{equation*}
\left|\frac{\pl^k}{\pl t^k}D^{\al}v\right|\le C_o 
\quad\text{ in }Q_\dl\text{ for }|\al|+2k\le4[\frac N2]+16,
\end{equation*}
where $[a]$ denotes the integer part of $a$. Thus in particular 
an upper bound on these derivatives up to the indicated order, 
gives their analyticity as signified by (\ref{Eq:An:7}). Assuming 
such an upper bound for the moment, we return to the 
original coordinates to get 
\begin{equation}\label{Eq:An:9}
\begin{aligned}
&|D^{\al}u\pto|=|D^{\al}v(0,0)|\frac{u\pto}{\rho^{|\al|}}
\le \frac{C H^{|\al|} |\al|!}{\rho^{|\al|}}\,u\pto,\\
&\left|\frac{\pl^k}{\pl t^k}u\pto\right|
=\left|\frac{\pl^k}{\pl t^k}v(0,0)\right|
\frac{u\pto^{1-k}}{\rho^{2k}}\le 
\frac{C H^{2k}(2k)!}{\rho^{2k}}\,u\pto^{1-k}.
\end{aligned}
\end{equation}
The proof is concluded, once the dependence of 
$C_o$ on $f(\eta)$ is determined. This estimation 
can be achieved by local DeGiorgi's or Moser's estimates. 
While the method is known, it is technically involved 
and reported in detail in Appendix~A. 

The analogous analyticity estimates for solutions to the 
porous medium type equation (\ref{Eq:4:1})--(\ref{Eq:4:2}) 
are similar, with the obvious changes, and we omit the 
details.
\appendix
\appendixpage
\section{Analyticity in the Space Variables, 
of Solutions to (\ref{Eq:1:1}). Estimating the first 
$4[\frac N2]+16$ Derivatives of $v$.}\label{App:An}
We will use expressions such as $w^4 f(w)$, $w^5 f^\prime$
and similar ones,
but we only have at our disposal the notion of weak solution, and
therefore, such a way of working does not seem justified. 
However, by the Harnack estimate of Theorem~\ref{Thm:3:1}, solutions are
classical, and in these calculations we are turning the \emph{qualitative}
information of $u$ being classical into the \emph{quantitative} information
of $u$ being analytic.

With respect to the previous sections, we use a different notation for cylinders, and we let 
$Q(\rho,\theta)=K_{\rho}\times(-\theta,0]$.
\subsection{An estimate of $\|Dv\|_{\infty}$}
Take $D_{x_i}$ of the logarithmic diffusion equation and set 
$w_i=v_{x_i}$ to get
\begin{equation*}
\pl_t w_{i}-\dvg\big(\frac1v Dw_i-\frac1{v^2}w_i Dv\big)=0
\end{equation*}
Setting $w=(w_1,\dots,w_N)$,
yields 
\begin{equation}\label{Eq:A:1}
w_t-\dvg\big(\frac1v Dw-\frac1{v^2}w\otimes w\big)=0.
\end{equation} 
For all derivations below we stipulate that
$\lm<1$, $\Lm=\lm^{-1}>1$,
$\lm<\frac1v<\Lm$, $\theta<1$ and $\rho<1$.
\begin{proposition}\label{Prop:A:1}
Let $w$ be a solution to (\ref{Eq:A:1}) and $\z$ be
a cutoff function in $Q=Q(\rho,\theta)$. Then
\begin{equation*}\label{Eq:A:2}
\begin{aligned}
&\sup_{-\theta<t<0}\int_{K_{\rho}}
\int_0^{|w|}sf(s)ds\z^2\,dx\\
&\quad +\frac{\lm}2\iint_Q |Dw|^2f(|w|)\z^2\,dxdt
+\frac{\lm}2\iint_Q|w||D|w||^2f'(|w|)\z^2\,dxdt\\
&\le\frac{2\Lm^2}{\lm}\iint_Q f(|w|)|w|^2|D\z|^2\,dxdt
+\frac{2\Lm^4}{\lm}\iint_Q |w|^4f(|w|)\z^2\,dxdt\\
&+\frac{\Lm^4}{2\lm}\iint_Q |w|^5f'(|w|)\z^2\,dxdt
+\iint_{Q}\int_0^{|w|}sf(s)ds2\z\z_t\,dxdt,
\end{aligned}
\end{equation*}
where $f:\rr_+\to\rr_+$ is a bounded, non-negative, 
non-decreasing Lipschitz continuous function.
\end{proposition}
\noi{\it Proof -} Multiply (\ref{Eq:A:1}) by $w f(|w|)\z^2$ and integrate by parts to get
\begin{equation*}
0=\iint_Q\frac{\pl}{\pl t}\int_0^{|w|}sf(s)ds\z^2\,dxdt+I
\end{equation*}
where
\begin{equation*}
\begin{aligned}
&I=\\
&=\iint_Q\left(\frac1v Dw-\frac1{v^2}w\otimes w\right)
\cdot\left[Dw f(|w|)\z^2+w f^\prime(|w|)D|w|\z^2\right.\\
&\left.+w f(|w|)2\z D\z\right]dxdt\\
&=\iint_Q\left[\frac1v f(|w|)|Dw|^2\z^2
+\frac1v w f^\prime(|w|)Dw\cdot D|w|\z^2
+\frac1v wf(|w|)2\z Dw\cdot D\z\right.\\
&\left.-\frac1{v^2}w\otimes w\cdot Dw f(|w|)\z^2
-\frac1{v^2}w\otimes w\cdot wf^\prime(|w|)D|w|\z^2\right.\\
&\left.-\frac1{v^2}w\otimes w\cdot wf(|w|)2\z D\z\right]dxdt\\
&\ge\iint_Q\left[\frac1v|Dw|^2f(|w|)\z^2
+\frac1v|w||D|w||^2f^\prime(|w|)\z^2
+\frac1v|w|f(|w|)2\z D|w|\cdot D\z\right.\\
&\left.-\frac1{v^2}|w|wD|w|f(|w|)\z^2
-\frac1{v^2}|w|^2wD|w|f^\prime(|w|)\z^2
-\frac1{v^2}|w|^2wD\z f(|w|)2\z\right]dxdt
\end{aligned}
\end{equation*}
Observing that $|D|w||\le |Dw|$ the previous inequality yields
\begin{equation*}
\begin{aligned}
&\iint_Q\left[\frac{\pl}{\pl t}\int_0^{|w|}sf(s)ds\z^2
+\lm|Dw|^2f(|w|)\z^2+\lm|w||D|w||^2f'(|w|)\z^2\right]dxdt\\
&\le \iint_Q\left[\Lm|w|f(|w|)2\z D|w|\cdot D\z
+\Lm^2|w|^2D|w|f(|w|)\z^2\right.\\
&+\left.\Lm^2|w|^3D|w|f'(|w|)\z^2
+\Lm^2|w|^3|D\z| f(|w|)2\z\right]dxdt\\
&\le \frac{\lm}2\iint_Q |Dw|^2f(|w|)\z^2\,dxdt 
+\frac{\Lm^2}{\lm}\iint_Q f(|w|)|w|^2|D\z|^2\,dxdt\\
&+\frac{\Lm^4}{\lm}\iint_Q |w|^4f(|w|)\z^2\,dxdt
+\Lm^2\iint_Q |w|^3|D\z| f(|w|)2\z\,dxdt\\
&+\frac{\lm}2\iint_Q |w||D|w||^2f^\prime(|w|)\z^2\,dxdt
+\frac{\Lm^4}{2\lm}\iint_Q |w|^5f^\prime(|w|)\z^2\,dxdt.\qquad\qquad\bbox
\end{aligned}
\end{equation*}
We will use this energy estimate and Moser's 
iteration to derive a bound for $\|Dv\|_{\infty}$.
\begin{proposition}
There exists a positive parameter $\dl$ that depends only 
on $\frac\Lm\lm$, such that if $\rho\in(0,\dl]$ then
\begin{equation*}
\|w\|_{\infty,Q(\sig\rho,\sig\theta)}\le 
\gm\bigg(\frac{\Lm}{\lm}\bigg)^{\mu_1}
\frac{(1+\theta^{-\mu_2})}{(1-\sig)^{\mu_2}}.
\end{equation*}
where $\gm$, $\mu_1$, and $\mu_2$ are positive parameters 
that depend only on $N$.
\end{proposition}
\noi{\it Proof -} We first estimate $\|w\|_2$. 
In the weak formulation of (\ref{Eq:1:1}) take 
the test function $v\z^2$. By standard calculations 
\begin{equation*}
\lm\iint_{Q(\rho,\theta)} |Dv|^2\z^2\,dxdt
\le \gm\int_{K_{\rho}\times\{-\theta\}}v^2\z^2\,dx
+\gm\iint_{Q(\rho,\theta)}v| D\z|^2\,dxdt.
\end{equation*}
Hence
\begin{equation}\label{Eq:A:3}
\iint_{Q(\sig\rho,\sig\theta)}|Dv|^2\,dxdt
\le \gm\frac{\Lm^2}{\lm}\frac{\rho^N}
{(1-\sig)^2}\left[1+\frac{\theta}{\rho^2}\right]i
\le\gm\left(\frac{\Lm}{\lm}\right)^2
\frac{1+\theta}{(1-\sig)^2\rho^2},
\end{equation}
The energy estimate \eqref{Eq:A:2} with 
$f(|w|)=|w|^{2\beta}$ and $\be\ge0$, yields
\begin{equation}\label{Eq:A:4}
\begin{aligned}
&\sup_{-\theta<t<0}\int_{K_{\rho}}|w|^{2\beta+2}\z^2(x,t)\,dx
+\iint_{Q(\rho,\theta)}|Dw|^2|w|^{2\beta}\z^2\,dxdt\\
&\le\frac{\Lm^4}{\lm^2}\iint_{Q(\rho,\theta)}
|w|^{2\beta+2}[|D\z|^2+\z_t]\,dxdt
+\frac{\Lm^4}{\lm^2}\iint_{Q(\rho,\theta)}
|w|^{2\beta+4}\z^2\,dxdt.
\end{aligned}
\end{equation}
Notice that $v$ is locally H\"older continuous with 
the H\"older norm and its exponent $\al$
depending only on $\frac{\Lm}{\lm}$ and $N$, that is
\begin{equation*}
|v(x,\cdot)-v(0,\cdot)|\le\gm\Lm\rho^\al.
\end{equation*}
Now we apply an integration by parts and use the 
H\"older continuity of $v$ to estimate the last 
term of \eqref{Eq:A:4}. 
If we freeze the time variable, then
\begin{equation}\label{Eq:A:5}
\begin{aligned}
I&=\int_{K_\rho}|w|^{2\be+4}\z^2\,dx
=\int_{K_\rho}|w|^{2\be+2}w\cdot w\z^2\,dx\\
&=\int_{K_{\rho}}|w|^{2\beta+2}\z^2Dv\cdot 
D(v-v(0,\cdot))\,dx\\
&=\int_{K_{\rho}}(v-v(0,\cdot))\left[\Delta v
|w|^{2\beta+2}\z^2+2|w|^{2\beta+2}\z D\z\cdot Dv\right.\\
&\left.+(2\beta+2)|w|^{2\beta+1}D|w|Dv\z^2\right]dx\\
&\le\gm\Lm\int_{K_{\rho}}\rho^{\al}
\left[\Delta v|w|^{2\beta+2}\z^2
+2|w|^{2\beta+2}\z D\z\cdot Dv\right.\\
&\left.+(2\beta+2)|w|^{2\beta+1}D|w|
\cdot Dv\z^2\right]dx\\
&\le \frac34\int_{K_{\rho}}|w|^{2\beta+4}\z^2\,dx
+\gm\rho^{2\al}\Lm^2\int_{K_{\rho}}|w|^{2\beta}|Dw|^2\,dx\\
&+\gm\rho^{2\al}\Lm^2\int_{K_{\rho}}|w|^{2\beta+2}|D\z|^2\,dx.
\end{aligned}
\end{equation}
Since all constants are independent of $t$, combining 
\eqref{Eq:A:4}--\eqref{Eq:A:5} yields
\begin{equation*}
\begin{aligned}
&\sup_{-\theta<t<0}\int_{K_{\rho}}|w|^{2\beta+2}\z^2(x,t)\,dx
+\iint_{Q(\rho,\theta)}|Dw|^2|w|^{2\beta}\z^2\,dxdt\\
&\le\frac{\Lm^4}{\lm^2}\iint_{Q(\rho,\theta)}|w|^{2\beta+2}[|D\z|^2+\z_t]dxdt\\
&+\gm\frac{\Lm^6}{\lm^2}\iint_{Q(\rho,\theta)}\rho^{2\al}\left[|w|^{2\beta}|Dw|^2+|w|^{2\beta+2}|D\z|^2\right]dxdt.
\end{aligned}
\end{equation*}
Thus, by taking
\begin{equation}\label{Eq:A:6}
\gm\frac{\Lm^6}{\lm^2}\rho^{2\al}=\frac12
\end{equation}
we have
\begin{equation*}
\begin{aligned}
&\sup_{-\theta<t<0}\int_{K_{\rho}}|w|^{2\beta+2}\z^2(x,t)\,dx
+\iint_{Q(\rho,\theta)}|Dw|^2|w|^{2\beta}\z^2\,dxdt\\
&\le\gm\frac{\Lm^6}{\lm^2}\iint_{Q(\rho,\theta)}
|w|^{2\beta+2}\left[|D\z|^2+\z_t
+\rho^{2\al}|D\z|^2\right]dxdt\\
&\le\gm\frac{\Lm^6}{\lm^2}\iint_{Q(\rho,\theta)}
|w|^{2\beta+2}\left[|D\z|^2+\z_t\right]dxdt,
\end{aligned}
\end{equation*}
as $\rho^{2\al}<1$ by \eqref{Eq:A:6}.
An application of the Sobolev parabolic embeddings gives
\begin{equation*}
\iint_{Q(\sig\rho,\sig\theta)}|w|^{(2\beta+2)\frac{N+2}N}dxdt\,
\le\left\{\gm\frac{\Lm^6}{\lm^2}
\iint_{Q(\rho,\theta)}|w|^{2\beta+2}\left[|D\z|^2+\z_t
\right]dxdt\right\}^{1+\frac2N}.
\end{equation*}
Now take 
\begin{equation*}
\rho_n=\sig\rho+\frac{1-\sig}{2^n}\rho,\qquad
\theta_n=\sig\theta+\frac{1-\sig}{2^n}\theta,\qquad
Q_n=Q(\rho_n,\theta_n), 
\end{equation*}
let $\z$ be a standard cutoff function in $Q_n$, and set 
\begin{equation*}
a_n=2\Big(\frac{N+2}N\Big)^n\quad \text{ and }\quad 
I_n=\iint_{Q_n}|w|^{a_n}dxdt.
\end{equation*}
Let $b=4^a$ and $a=1+\frac2N$. Begin from $2\be+2=a_o$, 
that is from $\be=0$, and apply the above estimate recursively 
up to $2\be+2=a_n$. This gives
\begin{equation*}
I_{n+1}\le C^ab^nI_n^{a},\quad\text{ where }\quad 
C=\frac{\Lm^6}{\lm^6}
\left\{\frac1{(1-\sig)^2\rho^2}+\frac1{(1-\sig)\theta}\right\}.
\end{equation*}
Iterating these recursive inequalities yields
\begin{equation*}
\begin{aligned}
I_{n+1}&\le C^{\sum_{i=0}^n a^{i+1}}b^{\sum_{i=0}^n(n-i)a^{i}}}
I_o^{a^{n+1}\\
&\le C^{\gm(N)a^{n+1}}b^{\gm(N)a^{n+1}}I_o^{a^{n+1}}.
\end{aligned}
\end{equation*}
Now take the  $\frac1{a^{n+1}}$ power of both sides 
and let $n\to\infty$. Taking into account the estimate of $I_o$ in 
\eqref{Eq:A:3}, we obtain
\begin{equation*}
\begin{aligned}
\|w\|_{\infty,Q(\sig\rho,\sig\theta)}&\le \gm
\bigg(\frac{\Lm}{\lm}\bigg)^{\mu_1}
\frac{(1+\theta^{-\mu_2})}{(1-\sig)^{\mu_2}}I_o\\
&\le \gm\bigg(\frac{\Lm}{\lm}\bigg)^{\mu_1}
\frac{(1+\theta^{-\mu_2})}{(1-\sig)^{\mu_2}},
\end{aligned}
\end{equation*}
where the constants $\mu_1$ and $\mu_2$ have been properly modified.
\hfill\bbox
\subsection{An Upper Bound for $\|v_t\|_{\infty}$}
\begin{proposition}
Let $v$ be a classical solution to the logarithmic diffusion equation
and assume $0<\lm\le v^{-1}\le\Lm$ in $Q(\rho,\theta)$; then
\begin{equation*}
\|v_t\|_{\infty,Q(\sig\rho,\sig\theta)}\le \gm\bigg(\frac{\Lm}{\lm}\bigg)^{\mu_1}\frac{(1+\theta^{-\mu_2})}{(1-\sig)^{\mu_2}},
\end{equation*}
where $\gm$, $\mu_1$, and $\mu_2$ are positive parameters that depend only on $N$.
\end{proposition}
\noi{\it Proof -} Multiply \eqref{Eq:1:1} by the test function $v_t\z^2$ and integrate over
the cylinder $Q_1=Q(\frac{\rho+\sig\rho}2,\frac{\theta+\sig\theta}{2})$, 
where $\sig\in(0,1)$.
Here $\z$ vanishes on the parabolic boundary of $Q_1$ and takes value $1$ in $Q(\sig\rho,\sig\theta)$.
A standard calculation gives
\begin{equation*}
\begin{aligned}
0&=\iint_{Q_1}\bigg[v^2_t\z^2dxdt+\frac1v Dv[\z^2Dv_t+2v_t\z D\z]\bigg]dxdt\\
&=\iint_{Q_1}\bigg[v^2_t\z^2dxdt+\frac1{2v}\z^2\frac{\pl}{\pl t}|Dv|^2+\frac2v v_t\z Dv D\z\bigg]dxdt\\
&=\iint_{Q_1}v^2_t\z^2 dxdt+\int_{K_{\frac{\rho+\sig\rho}2}\times \{0\}}\frac1{2v}|Dv|^2\z^2\,dx\\
&-\iint_{Q_1}|Dv|^2\bigg[-\frac1{2v^2}v_t\z^2+\frac1v\z \z_t\bigg]dxdt+\iint_{Q_1}\frac2v v_t\z Dv D\z dxdt.
\end{aligned}
\end{equation*}
This gives the estimate
\begin{equation*}
\begin{aligned}
\iint_{Q_1}v_t^2\z^2dxdt&\le \iint_{Q_1}|Dv|^2\frac1{2v^2}|v_t|\z^2dxdt+\iint_{Q_1}\frac1v\z|\z_t||Dv|^2dxdt\\
&+\iint_{Q_1}\frac1v|Dv||v_t|2\z|D\z|dxdt\\
&\le \frac12\iint_{Q_1}v_t^2\z^2dxdt+\iint_{Q_1}\frac{|Dv|^4}{4v^4}\z^2dxdt\\
&+\iint_{Q_1}\frac{4|Dv|^2}{v^2}|D\z|^2dxdt+\iint_{Q_1}\frac{|Dv|^2}{v}|\z_t|dxdt.
\end{aligned}
\end{equation*}
Taking into account the estimate for $\|Dv\|_{\infty,Q_1}$ of the previous section,
we have 
\begin{equation*}
\|v_t\|_{2,Q(\sig\rho,\sig\theta)}\le \gm\bigg(\frac{\Lm}{\lm}\bigg)^{\mu_1}\frac{(1+\theta^{-\mu_2})}{(1-\sig)^{\mu_2}}.
\end{equation*}
for some $\mu_1(N)$, and $\mu_2(N)>0$.

Now take the time derivative of the logarithmic diffusion equation, 
and in the corresponding weak formulation use the test function $v_tf(|v_t|)\z^2$
where $f:\rr_+\to\rr_+$ is a bounded, non-decreasing Lipschitz function,
and $\z$ vanishes on the parabolic boundary of $Q=Q(\rho,\theta)$
and takes value $1$ in $Q(\sig\rho,\sig\theta)$.
Let $M=\|Dv\|_{\infty,Q}$. A standard calculation yields
\begin{equation*}
\begin{aligned}
&\iint_{Q}\frac12\frac{\pl}{\pl t}|v_t|^2f\z^2dxdt+\lm\iint_{Q}|Dv_t|^2f\z^2dxdt\\
&+\lm\iint_{Q}|Dv_t|^2|v_t|f'\z^2dxdt\\
&\le 2\Lm\iint_{Q}|Dv_t||v_t|f\z|D\z|dxdt\\
&+\Lm^2\iint_{Q}|v_t| |Dv|\left[|Dv_t|f\z^2
+|v_t|f'|D v_t|\z^2+2|v_t|f\z |D\z|\right]dxdt\\
&=I_1+I_2+I_3+I_4.
\end{aligned}
\end{equation*}
Let us estimate the four terms.
\begin{equation*}
\begin{aligned}
I_1&\le \frac{\lm}4\iint_{Q}|Dv_t|^2f\z^2dxdt+\frac{16\Lm^2}{\lm}\iint_{Q}|v_t|^2f|D\z|^2dxdt;\\
I_2&\le \frac{\lm}4\iint_{Q}|Dv_t|^2f\z^2dxdt+4\frac{M^2\Lm^4}\lm\iint_{Q}|v_t|^2f \z^2dxdt;\\
I_3&\le \frac{\lm}2\iint_{Q}|Dv_t|^2 |v_t|f'\z^2dxdt+4\frac{M^2\Lm^4}{\lm}\iint_{Q}|v_t|^3f' \z^2dxdt;\\
I_4&\le M\Lm^2\iint_{Q}|v_t|^2f \z^2dxdt+M\Lm^2\iint_{Q}|v_t|^2 f |D\z|^2dxdt.
\end{aligned}
\end{equation*}
Summarizing we have 
\begin{equation*}
\begin{aligned}
&\sup_{-\theta<t<0}\int_{K_{\rho}}\int_0^{|v_t|}sf(s)ds\z^2dx+\lm\iint_Q|Dv_t|^2f\z^2dxdt\\
&+\lm\iint_Q|D|v_t||^2f'|v_t|\z dxdt\\
&\le\gm\bigg(\frac{\Lm}{\lm}\bigg)^4[M^2+1]\iint_Q\left[|v_t|^2f\z^2+|v_t|^3f'\z^2+|v_t|^2 f|D\z|^2\right]dxdt\\
&+2\iint_Q\int_0^{|v_t|}sf(s)ds\,\z\z_tdxdt.
\end{aligned}
\end{equation*}
Now take $f(s)=s^\be$ for $\be\ge0$; then
\begin{equation*}
\begin{aligned}
&\frac1{\be+2}\sup_{-\theta<t<0}\int_{K_{\rho}}|v_t|^{\be+2}\z^2dx+\iint_Q|D|v_t||^2|v_t|^{\be}\z^2dxdt\\
&\le\gm\bigg(\frac{\Lm}{\lm}\bigg)^5[M^2+1]\iint_Q|v_t|^{\be+2}\left[1+\be+|D\z|^2+|\z_t|\right]dxdt.
\end{aligned}
\end{equation*}

\noi Let $a=1+\frac2N$, 
$$C=\gm\bigg(\frac{\Lm}{\lm}\bigg)^5[M^2+1]\left[1+\frac1{(1-\sig)^2\rho^2}+\frac1{(1-\sig)\theta}\right],$$
 and let $w\df=v_t$.
An application of the Sobolev embedding yields
\begin{equation*}
\begin{aligned}
&\iint_{Q(\sig\rho,\sig\theta)}|w|^{(\be+2)\frac{N+2}N}dxdt\\
&\le \bigg(\sup_{-\theta<t<0}\int_{K_{\rho}}||w|^{\frac{\be+2}2}\z|^2dx\bigg)^{\frac2N}
\iint_{Q}|D[|w|^{\frac{\be+2}2}\z]|^2dxdt\\
&\le C^a (1+\be)^{3a}\bigg(\iint_Q|w|^{\be+2}dxdt\bigg)^a.
\end{aligned}
\end{equation*}
Take $Q_n$ as before and define
$$\be_o=0,\ \ \be_{n+1}+2=(\be_n+2)\frac{N+2}N\ \ \Rightarrow\ \  \be_n=2\bigg(\frac{N+2}N\bigg)^n-2,$$
and
$$ I_n=\iint_{Q_n}|w|^{\be_n+2}dxdt.$$
It then follows that
\begin{equation*}
I_{n+1}\le C^ab^nI_n^a
\end{equation*}
for some positive constant $b$ depending only on $N$.
A standard iteration  gives
\begin{equation*}
\begin{aligned}
I_{n+1}&\le C^{\sum_{i=1}^{n+1}a^i}b^{\sum_{i=0}^{n}(n-i)a^i}I_o^{a^{n+1}}\\
&\le C^{\gm(N) a^{n+1}}b^{\gm(N) a^{n+1}}I_o^{a^{n+1}}.
\end{aligned}
\end{equation*}
Therefore, taking the  $\frac1{a^{n+1}}$ power of both sides and letting $n\to\infty$, we have
\begin{equation*}
\|v_t\|_{\infty, Q(\sig\rho,\sig\theta)}\le C^{\gm}b^{\gm}\iint_Q|v_t|^2dxdt.
\end{equation*}

\noi To conclude this section, bound the right hand side using the estimate for $\|v_t\|_2$
on an intermediate cylinder.
\hfill\bbox
\subsection{An Upper Bound for $\frac{\pl^k}{\pl t^k}D^{\al}v$}
Differentiating  the logarithmic diffusion equation successively we have
\begin{equation}\label{Eq:A:7}
\begin{aligned}
&\frac{\pl^k}{\pl t^k}D^{\al}v_t-\dvg\bigg(\frac1vD\frac{\pl^k}{\pl t^k}D^{\al}v
+\sum_{j<k}\binom{k}{j}\frac{\pl^{k-j}}{\pl t^{k-j}}\frac1v\frac{\pl^j}{\pl t^j}DD^{\al}v\\
&+\sum_{|\beta|<|\al|}\binom{\al}{\be}\sum_{j\le k}
\binom{k}{j}\frac{\pl^{k-j}}{\pl t^{k-j}}D^{\al-\be}\frac1v\frac{\pl^j}{\pl t^j}D^{\be}Dv
\bigg)=0,
\end{aligned}
\end{equation}
where $k\in\nn$ and $\al$ is a multi-index.
\noi For an integer $n>0$ let 
\begin{equation*}
|w|^2=\sum_{k+|\al|=n}\left|\frac{\pl^k}{\pl t^k}D^{\al}v\right|^2.
\end{equation*}
We have the following bound for general derivatives of the logarithmic diffusion equation.
\begin{proposition}\label{Prop:A:4}
Let $v$ be a classical solution to the logarithmic diffusion equation in 
$Q(\rho,\theta)$ and fix $\sig\in(0,1)$. Assume
$0<\lm\le v^{-1}\le\Lm$ in $Q(\rho,\theta)$. There exists a 
positive parameter $\dl$ that depends only  $\frac\Lm\lm$, such that if 
$\rho\in(0,\dl]$, then in $Q(\sig\rho,\sig\theta)$ 
\begin{equation}\label{Eq:A:11}
\|w\|_{\infty,Q(\sig\rho,\sig\theta)}\le \gm\bigg(\frac{\Lm}{\lm}\bigg)^{\mu_1}\frac{(1+\theta^{-\mu_2})}{(1-\sig)^{\mu_2}},
\end{equation}
where $\gm$, $\mu_1$ and $\mu_2$ are positive parameters that depend only on $N$ and $n$.
\end{proposition}
\noi{\it Proof -} Multiply \eqref{Eq:A:7} by the test function $\frac{\pl^k}{\pl t^k}D^{\al}vf(|w|)\z^2$,
where $f:\rr_+\to\rr_+$ is a bounded, non decreasing Lipschitz function.
Here $\z$ vanishes on the parabolic boundary of 
$Q_2=Q({\frac{(1+\sig)\rho}2},\frac{(1+\sig)\theta}2)$ and takes value $1$ 
in $Q(\sig\rho,\sig\theta)$.
Standard calculations and a sum over $k+|\al|=n$ give
\begin{equation*}
\begin{aligned}
&\iint_{Q_2}\z^2\frac{\pl}{\pl t}\int_0^{|w|}sf(s)ds\,dxdt+\lm\iint_{Q_2}|w||D|w||^2f'(|w|)\z^2dxdt\\
&+\lm\sum_{k+|\al|=n}\iint_{Q_2}|D\frac{\pl^k}{\pl t^k}D^{\al}v|^2f(|w|)\z^2dxdt
\le I,
\end{aligned}
\end{equation*}
where
\begin{equation*}
\begin{aligned}
I&=-\sum_{k+|\al|=n}\iint_{Q_2}\frac2v \z f(|w|) \frac{\pl^k}{\pl t^k}D^{\al}v\, 
D\frac{\pl^k}{\pl t^k}D^{\al}v  D\z dxdt\\
&-\sum_{k+|\al|=n}\iint_{Q_2}\bigg[\sum_{j<k}\binom{k}{j}\frac{\pl^{k-j}}{\pl t^{k-j}}\frac1v\frac{\pl^j}{\pl t^j}DD^{\al}v\\
&+\sum_{|\beta|<|\al|}\binom{\al}{\be}\sum_{j}
\binom{k}{j}\frac{\pl^{k-j}}{\pl t^{k-j}}D^{\al-\be}\frac1v\frac{\pl^j}{\pl t^j}D^{\be}Dv\bigg]\\
&\times\bigg[f(|w|) \z^2 D\frac{\pl^k}{\pl t^k}D^{\al}v +\z^2\frac{\pl^k}{\pl t^k}D^{\al}v f'(|w|) D|w|+
2\z f(|w|) \frac{\pl^k}{\pl t^k}D^{\al}v D\z\bigg]dxdt.
\end{aligned}
\end{equation*}
Notice that
\begin{equation*}
\begin{aligned}
&\left|\sum_{j<k}\binom{k}{j}\frac{\pl^{k-j}}{\pl t^{k-j}}\frac1v\frac{\pl^j}{\pl t^j}DD^{\al}v\right|\\
&+\left|\sum_{|\beta|<|\al|}\binom{\al}{\be}\sum_{j\le k}
\binom{k}{j}\frac{\pl^{k-j}}{\pl t^{k-j}}D^{\al-\be}\frac1v\frac{\pl^j}{\pl t^j}D^{\be}Dv\right|\\
&\le P[|w|+1]
\end{aligned}
\end{equation*}
for some polynomial $P$ with variables $\{\|\frac{\pl^k}{\pl t^k}D^{\al}v\|_{\infty,Q_2} \text{ for }k+|\al|<n;\frac{\Lm}{\lm}\}$.
Thus 
\begin{equation*}
\begin{aligned}
I&\le \frac{\lm}4\sum_{k+|\al|=n}\iint_{Q_2}|D\frac{\pl^k}{\pl t^k}D^{\al}v|^2f(|w|)\z^2dxdt\\
&+\frac{\Lm^2}{2\lm}\sum_{k+|\al|=n}\iint_{Q_2}|\frac{\pl^k}{\pl t^k}D^{\al}v|^2 f(|w|) |D\z|^2dxdt\\
&+\frac{\lm}4\sum_{k+|\al|=n}\iint_{Q_2}|D\frac{\pl^k}{\pl t^k}D^{\al}v|^2 f(|w|) \z^2dxdt\\
&+\frac\gm\lm
\iint_{Q_2}P^2(|w|^2+1)f(|w|) \z^2dxdt\\
&+\frac{\lm}2\iint_{Q_2}|D|w||^2f'(|w|) |w|\z^2dxdt+\frac\gm\lm\iint_{Q_2}P^2(|w|^2+1)|w|f'(|w|) \z^2dxdt\\
&+\iint_{Q_2}P(|w|+1)|w| f(|w|) \z |D\z| dxdt.
\end{aligned}
\end{equation*}
We obtain the following energy estimate
\begin{equation*}
\begin{aligned}
&\sup_{-\frac{\theta+\sig\theta}2<t<0}\int_{K_{\frac{\rho+\sig\rho}2}}\int_0^{|w|}sf(s)ds\z^2dx+
\frac\lm2\iint_{Q_2}|w||D|w||^2f'(|w|)\z^2dxdt\\
&+\frac\lm2\sum_{k+|\al|=n}\iint_{Q_2}|D\frac{\pl^k}{\pl t^k}D^{\al}v|^2f(|w|)\z^2dxdt\\
&\le\gm\left(\frac\Lm\lm\right)^2\left[\iint_{Q_2}P^2(|w|^2+1)f(|w|) \z^2dxdt+\iint_{Q_2}|w|^2 f(|w|) |D\z|^2dxdt\right.\\
&\left.+\iint_{Q_2}P^2(|w|^2+1)|w| f'(|w|) \z^2 dxdt\right]+2\iint_{Q_2}\int_0^{|w|}sf(s)ds\,\z|\z_t|dxdt,\\
\end{aligned}
\end{equation*}
and also
\begin{equation*}
\begin{aligned}
&\sup_{-\frac{\theta+\sig\theta}2<t<0}\int_{K_{\frac{\rho+\sig\rho}2}}\int_0^{|w|}sf(s)ds\z^2dx+
{\frac\lm2}\iint_{Q_2}{|w|}|D|w||^2f'(|w|)\z^2dxdt\\
&+{\frac\lm2}\iint_{Q_2}|D|w||^2f(|w|)\z^2dxdt\\
&\le{\gm\left(\frac\Lm\lm\right)^2}\left[\iint_{Q_2}P^2(|w|^2+1)f(|w|) \z^2dxdt+\iint_{Q_2}|w|^2 f(|w|) |D\z|^2dxdt\right.\\
&{\left.+\iint_{Q_2}P^2(|w|^2+1)|w| f'(|w|) \z^2 dxdt\right]}+2\iint_{Q_2}\int_0^{|w|}sf(s)ds\,\z|\z_t|dxdt.
\end{aligned}
\end{equation*}
Now we have at our disposal the  sup-estimates for $Dv$ and $v_t$ in terms of $\Lm/\lm$ only.
Next we assume the supremum of $\frac{\pl^k}{\pl t^k}D^{\al}v$ is estimated
for all $k+|\al|<n$ by a similar quantity as the right hand side of \eqref{Eq:A:11}, 
By Moser's method, the above energy estimate will yield
a bound for the
case $k+|\al|=n$. These will depend on the $L^2$ norms of $\frac{\pl^k}{\pl t^k}D^{\al}v$
and a polynomial with variables $\{\|\frac{\pl^k}{\pl t^k}D^{\al}v\|_{\infty,Q_2} \text{ for }k+|\al|<n;\frac{\Lm}{\lm}\}$.

Take $f(s)=s^{\be}$ for $\be\ge0$; then the energy estimate yields
\begin{equation*}
\begin{aligned}
&\frac1{\be+2}\sup_{-\frac{\theta+\sig\theta}2<t<0}\int_{K_{\frac{\rho+\sig\rho}2}}|w|^{\be+2}\z^2dx+\iint_{Q_2}|D|w||^2|w|^{\be}\z^2dxdt\\
&\le \gm \left(\frac\Lm\lm\right)^3P^2(1+\be)\bigg[1+\frac1{(1-\sig)^2\rho^2}+\frac1{(1-\sig)\theta}\bigg]\iint_{Q_2}[|w|^{\be+2}+|w|^{\be}]dxdt.
\end{aligned}
\end{equation*}
Let $a=1+\frac2N$ and $C=\gm \left(\frac\Lm\lm\right)^3\bigg[1+\frac1{(1-\sig)^2\rho^2}+\frac1{(1-\sig)\theta}\bigg]$
and assume $|Q_2|<1$;
 then an application of the Sobolev embedding yields
\begin{equation*}
\begin{aligned}
&\iint_{Q(\sig\rho,\sig\theta)}|w|^{(\be+2)\frac{N+2}N}dxdt\\
&\le \bigg(\sup_{-\frac{\theta+\sig\theta}2<t<0}\int_{K_{\frac{\rho+\sig\rho}2}}||w|^{\frac{\be+2}2}\z|^2dx\bigg)^{\frac2N}
\iint_{Q_2}|D[|w|^{\frac{\be+2}2}\z]|^2dxdt\\
&\le C^aP^{2a}(1+\be)^{3a}\bigg[\iint_{Q_2}(|w|^{\be+2}+1)dxdt\bigg]^a\\
&\le C^aP^{2a} (1+\be)^{3a}\bigg(\iint_{Q_2}|w|^{\be+2}dxdt\bigg)^a+C^aP^{2a} (1+\be)^{3a},
\end{aligned}
\end{equation*}
after a proper adjustment of the constant $\gm$ in the definition of $C$.
Take 
\begin{equation*}
\rho_n=\sig\rho+\frac{1-\sig}{2^{n+1}}\rho,\qquad\theta_n=\sig\theta+\frac{1-\sig}{2^{n+1}}\theta,
\qquad Q_n=Q({\rho_n},\theta_n).
\end{equation*}
Define
$$\be_o=0,\ \ \be_{n+1}+2=(\be_n+2)\frac{N+2}N\ \ \Rightarrow\ \  \be_n=2\bigg(\frac{N+2}N\bigg)^n-2,$$
and
$$ I_n=\iint_{Q_n}|w|^{\be_n+2}dxdt.$$
We have that
\begin{equation*}
I_{n+1}\le C^aP^{2a}b^nI_n^a+C^aP^{2a}b^n,
\end{equation*}
where the constant $b$ depends only on $N$.
A standard iteration and a proper adjustment of $P$ give
\begin{equation*}
\begin{aligned}
I_{n+1}&\le 2^{\sum_{i=1}^n(a^i-1)}(CP^{2})^{\sum_{i=1}^{n+1}a^i}b^{\sum_{i=0}^n(n-i)a^i}I_o^{a^{n+1}}\\
&+2^{\sum_{i=1}^n(a^i-1)}(CP^{2})^{\sum_{i=1}^{n+1}a^i}b^{\sum_{i=0}^n(n-i)a^i}\\
&\le (CP^{2})^{\gm(N) a^{n+1}}I_o^{\gm(N) a^{n+1}}+(CP^{2})^{\gm(N) a^{n+1}}
\end{aligned}
\end{equation*}
Take  the $\frac1{a^{n+1}}$  power of both sides and let $n\to\infty$ to obtain
\begin{equation}\label{Eq:A:10}
\|w\|_{\infty, Q(\sig\rho,\sig\theta)}\le (CP^{2})^{\gm}\iint_{Q_2}|w|^2dxdt+(CP^{2})^{\gm}.
\end{equation}
Remember that $Q_2=Q(\frac{\rho+\sig\rho}2,\frac{\theta+\sig\theta}2)$; in order 
to conclude the proof of Proposition~\ref{Prop:A:4}, we only need 
to estimate $\|w\|_{2,Q_2}$.
\subsubsection{An Estimate of $\|w\|_{2,Q_2}$}
It is enough to give an estimate 
of $\|w\|_{2,Q(\sig\rho,\sig\theta)}$.
Replace $k$ in \eqref{Eq:A:7} by $k-1$ and assume $k\ge1$.
We can rewrite \eqref{Eq:A:7} as 
\begin{equation*}
\frac{\pl^{k-1}}{\pl t^{k-1}}D^{\al}v_t-\dvg\bigg(\frac1vD\frac{\pl^{k-1}}{\pl t^{k-1}}D^{\al}v\bigg)=\dvg f
\end{equation*}
where  
\begin{align*}
f&=\sum_{j<k-1}\binom{k-1}{j}\frac{\pl^{k-1-j}}{\pl t^{k-1-j}}\frac1v\frac{\pl^j}{\pl t^j}DD^{\al}v\\
&+\sum_{|\beta|<|\al|}\binom{\al}{\be}\sum_{j\le k-1}
\binom{k-1}{j}\frac{\pl^{k-1-j}}{\pl t^{k-1-j}}D^{\al-\be}\frac1v\frac{\pl^j}{\pl t^j}D^{\be}Dv
\end{align*}
If $\z$ is a smooth function in $Q_2$
and takes value $1$ in $Q(\sig\rho,\sig\theta)$, then the standard $L^2$ estimate for the linear parabolic equations gives
\begin{equation}\label{Eq:A:12}
\begin{aligned}
&\|\frac{\pl^k}{\pl t^k}D^{\al}v\|^2_{2,Q(\sig\rho,\sig\theta)}\\
&\le \gm \bigg[\|\dvg f\|^2_{2,Q_2}
+\bigg(\frac1{(1-\sig)\theta}+\frac1{(1-\sig)^2\rho^2}\bigg)\|\frac{\pl^{k-1}}{\pl t^{k-1}}D^{\al}v\|^2_{2,Q_2}\bigg].
\end{aligned}
\end{equation}
Let us denote $P$ as a polynomial of variables $\{\|\frac{\pl^k}{\pl t^k}D^{\al}v\|_{\infty,Q_2}\ \text{ for }k+|\al|<n;\frac{\Lm}{\lm}\}$.
Observe that
\begin{align*}
\|\dvg f\|_{2,Q_2}^2\le\gm P[I_1+I_2+I_3+1]
\end{align*}
where
\begin{align*}
&I_1=\iint_{Q_2}\left|\frac{\pl^{k-1}}{\pl t^{k-1}}Dv\right|^2\,dxdt;\\
&I_2=\iint_{Q_2}\left|\frac{\pl^{k-2}}{\pl t^{k-2}}D^{|\al|+2}v\right|^2\,dxdt;\\
&I_3=\iint_{Q_2}\left|\frac{\pl^{k-1}}{\pl t^{k-1}}D^{|\al|+1}v\right|^2\,dxdt.
\end{align*}
Here for an integer $l$
\[
|D^lv|^2=\sum_{|\be|=l}|D^{\be}v|^2.
\]
These quantities can all be estimated in the same way. Indeed, they all contain
spatial derivatives and we can use the principal part of the differentiated \eqref{Eq:1:1}
to estimate them. Precisely, we write \eqref{Eq:A:7} as
\begin{align*}
&\frac{\pl^{k-s}}{\pl t^{k-s}}D^{\eta}v_t-\dvg\bigg(\frac1vD\frac{\pl^{k-s}}{\pl t^{k-s}}D^{\eta}v
+\sum_{j<k-s}\binom{k-s}{j}\frac{\pl^{k-s-j}}{\pl t^{k-s-j}}\frac1v\frac{\pl^j}{\pl t^j}DD^{\eta}v\\
&+\sum_{|\beta|<|\eta|}\binom{\eta}{\be}\sum_{j\le k-s}
\binom{k-s}{j}\frac{\pl^{k-s-j}}{\pl t^{k-s-j}}D^{\eta-\be}\frac1v\frac{\pl^j}{\pl t^j}D^{\be}Dv
\bigg)=0,
\end{align*}
where $1\le s\le k$ and $|\eta|=|\al|+s$. 
Take the test function $$\frac{\pl^{k-s}}{\pl t^{k-s}}D^{\eta}v\z^2$$
where $\z$ vanishes on the parabolic boundary of 
$Q_3=Q(\frac{3(1+\sig)\rho}4,\frac{3(1+\sig)\theta}4)$ 
and takes value $1$ in $Q_2$.
Integrating in $Q_3$, a standard calculation yields
\begin{align*}
&\iint_{Q_3}\frac1v|D\frac{\pl^{k-s}}{\pl t^{k-s}}D^{\eta}v|^2\z^2\,dxdt\\
&=-\frac12\iint_{Q_3}\frac{\pl}{\pl t}\left|\frac{\pl^{k-s}}{\pl t^{k-s}}D^{\eta}v\right|^2\z^2\,dxdt\\
&-\iint_{Q_3}\bigg[\sum_{j<k-s}\binom{k-s}{j}\frac{\pl^{k-s-j}}{\pl t^{k-s-j}}\frac1v\frac{\pl^j}{\pl t^j}DD^{\eta}v\\
&+\sum_{|\beta|<|\eta|}\binom{\eta}{\be}\sum_{j\le k-s}
\binom{k-s}{j}\frac{\pl^{k-s-j}}{\pl t^{k-s-j}}D^{\eta-\be}\frac1v\frac{\pl^j}{\pl t^j}D^{\be}Dv\bigg]\\
&\times\bigg[D\frac{\pl^{k-s}}{\pl t^{k-s}}D^{\eta}v\z^2+2\frac{\pl^{k-s}}{\pl t^{k-s}}D^{\eta}v\z D\z\bigg]\\
&\le \frac\lm{2}\iint_{Q_3}|D\frac{\pl^{k-s}}{\pl t^{k-s}}D^{\eta}v|^2\z^2\,dxdt+\gm\bigg[\frac1{(1-\sig)\theta}+\frac1{(1-\sig)^2\rho^2}\bigg]P.
\end{align*}
This together with \eqref{Eq:A:12} gives 
\begin{equation}\label{Eq:A:13}
\|\frac{\pl^k}{\pl t^k}D^{\al}v\|^2_{2,Q(\sig\rho,\sig\theta)}\le \gm\bigg[\frac1{(1-\sig)\theta}+\frac1{(1-\sig)^2\rho^2}\bigg]P.
\end{equation}
Now the only remaining case is $|\al|=n$. For this, we consider
the equation \eqref{Eq:A:7} with $k=0$ and assume $|\al|=n-1$.
Take the test function $\z^2 D^{\al}v$, where $\z$ vanishes on the parabolic boundary
of $Q_3$ and takes $1$ in $Q(\sig\rho,\sig\theta)$.
\begin{equation*}
\begin{aligned}
\iint_{Q_3}\frac1v|DD^{\al}v|^2\z^2&dxdt=-\iint_{Q_3}\frac12\frac{\pl}{\pl t}(D^{\al}v)^2\z^2dxdt\\
&-\iint_{Q_3}\frac1v DD^{\al}vD^{\al}v2\z D\z dxdt\\
&-\sum_{|\be|<|\al|}\iint_{Q_3}\binom{\al}{\be}D^{\al-\be}\frac1vDD^{\be}vDD^{\al}v\z^2 dxdt\\
&-\sum_{|\be|<|\al|}\iint_{Q_3}\binom{\al}{\be}D^{\al-\be}\frac1vDD^{\be}vDD^{\al}vD^{\al}2\z D\z dxdt\\
&\le\frac{\lm}2\iint_{Q_3}|DD^{\al}v|^2\z^2dxdt+\bigg[\frac1{(1-\sig)\theta}+\frac1{(1-\sig)^2\rho^2}\bigg]P
\end{aligned}
\end{equation*}
Summing over all $|\al|=n-1$ actually gives
$$\sum_{|\al|=n}\iint_{Q_3}|D^{\al}v|^2\z^2dxdt\le\bigg[\frac1{(1-\sig)\theta}+\frac1{(1-\sig)^2\rho^2}\bigg]P.$$
If we take into consideration an intermediate cylinder, then this, together with \eqref{Eq:A:13} in \eqref{Eq:A:10}, yields
\begin{equation}
\|w\|_{\infty,Q(\sig\rho,\sig\theta)}\le P^{\gm}\bigg[\frac1{(1-\sig)\theta}+\frac1{(1-\sig)\rho}\bigg]^\gm
\end{equation}  
for some $\gm$ depending only on $N$. The induction hypothesis and the
definition of $\rho$ in \eqref{Eq:A:6} imply that
\begin{equation}
\|w\|_{\infty,Q(\sig\rho,\sig\theta)}\le \gm_1(N,n)\bigg(\frac{\Lm}{\lm}\bigg)^{\mu_1}\frac{1+\theta^{-\mu_2}}{(1-\sig)^{\mu_2}}.
\end{equation}
\hfill \bbox
\section{Proof of Proposition~\ref{Prop:4:1} for Weak 
Solutions to  
Equations (\ref{Eq:4:5})--(\ref{Eq:4:6})}\label{App:PM}
\subsection*{An Auxiliary Lemma}\label{App:PM:1}
\begin{lemma}\label{Lm:PM:1} 
Let $u$ be a non-negative, local, weak solution to 
the singular equations (\ref{Eq:4:5})--(\ref{Eq:4:6}), 
in $E_T$.  There exist two positive constants 
$\gm_1$, $\gm_2$ depending only on the 
data $\data$, such that for all cylinders 
$K_{4\rho}(y)\times[s,t]\subset E_T$, and all 
$\sig\in(0,1)$, 
\begin{equation*}
\begin{aligned}
\int_s^t\int_{K_\rho(y)}\frac{|Du|^2}{u^{2-\frac m2}}\z^2dx\,d\tau
&\le\gm_1(1+\Lm_{\frac m2,1})\rho^{N\frac m2}\,
\mcl{S}_\sig^{1-\frac m2}+\frac{\gm_2}{\sig^2\rho^2}
(\Lm^2_{\frac m2,1}\\
&\quad +\Lm^2_{\frac m2,2})\,\mcl{S}_\sig^{\frac m2}\,
(t-s)\rho^{N(1-\frac m2)},
\end{aligned}
\end{equation*}
where 
\begin{equation*}
\mcl{S}_\sig =\sup_{s<\tau<t}\int_{K_{(1+\sig)\rho}(y)} 
u(\cdot,\tau)dx.
\end{equation*}
\end{lemma}
{\it Proof -} In the following we restrict to $0<m<\frac23$, 
since we are mainly interested in proving the stability 
of the estimates as $m\to0^+$. For $m\in(\frac13,1)$ 
similar arguments hold, provided a slightly different 
test function $\vp$ is chosen (see \cite{DBGV-mono}, 
\S~B.1.1 for more details).

Assume $(y,s)=(0,0)$, fix $\sig\in(0,1)$,
and let $x\to\z(x)$ be a non-negative piecewise smooth 
cutoff function in $K_{(1+\sig)\rho}$ that vanishes outside 
$K_{(1+\sig)\rho}$, equals one on $K_\rho$, and such that 
$|D\z|\le(\sig\rho)^{-1}$. 
Let $s_1\in[0,t]$ be such that
\begin{equation*}
\ssig =\sup_{0 <s<t}\int_{K_{(1+\sig)\rho}(y)} 
u(\cdot,s)dx=\int_{K_{(1+\sig)\rho}(y)} 
u(\cdot,s_1)dx,
\end{equation*}
and set
\begin{equation*}
\ssigbar \df=\frac{\mcl{S}_\sig}{\rho^N}.
\end{equation*}
In the weak formulation of (\ref{Eq:4:5})--(\ref{Eq:4:6}) 
take the test function 
\begin{equation*}
\vp=\left(\testfunlu\right)_+ \z^2,
\end{equation*}
and integrate 
over $Q=K_{(1+\sig)\rho}\times(0,t]$, to obtain
\begin{align*}
0=&\iint_{Q}\frac{\pl}{\pl\tau} u \left(\testfunlu\right)_+
\z^2dx\,d\tau\\
&\quad+\iint_{Q} \bl{A}(x,\tau,u,Du)\cdot 
D\left[\left(\testfunlu\right)_+\z^2\right] dx\,d\tau\\
&= I_1+I_2.
\end{align*}
We estimate these two terms separately. 
\begin{equation*}
\begin{aligned}
I_1&=\iint_{Q} \frac{\pl}{\pl\tau} u 
\left(\testfunlu\right)_+\z^2dx\,d\tau\\
&=\iint_{Q\cap[u<\ssigbar]} \frac{\pl}{\pl\tau} u 
\left(\testfunlu\right)\z^2dx\,d\tau\\
&=-\int_{K_{(1+\sig)\rho}\cap[u<\ssigbar]}\z^2(x)
\left(\int_{u(x,t)}^{\ssigbar}\testfunls\,ds\right)dx\\
&\quad +\int_{K_{(1+\sig)\rho}\cap[u<\ssigbar]}\z^2(x)
\left(\int_{u(x,0)}^{\ssigbar}\testfunls\,ds\right)dx.
\end{aligned}
\end{equation*}
Next,
\begin{align*} 
I_2&=\iint_{Q}\bl{A}(x,\tau,u,Du)\cdot D
\left[\left(\testfunlu\right)_+\z^2\right] dx\,d\tau\\
&=\iint_{Q\cap[u<\ssigbar]} \bl{A}(x,\tau,u,Du)\cdot 
D\left[\left(\testfunlu\right)\z^2\right] dx\,d\tau\\
&=-\frac12\iint_{Q\cap[u<\ssigbar]}\z^2 
u^{-\frac m2-1}\bl{A}(x,\tau,u,Du)\cdot Du\, dx\, d\tau\\
&\quad +2\iint_{Q\cap[u<\ssigbar]}\z\,
\left(\testfunlu\right)\bl{A}(x,\tau,u,Du)\cdot{D\z}\, dx\, d\tau\\
&\le-\frac{C_o}2\iint_{Q\cap[u<\ssigbar]}
u^{-\frac m2-1}u^{m-1}\z^2|Du|^2\,dx\,d\tau\\
&\quad+2C_1\iint_{Q\cap[u<\ssigbar]}\z
\left(\testfunlu\right) u^{m-1}|Du| |D\z| dx\, d\tau\\
&\le-\frac{C_o}4\iint_{Q\cap[u<\ssigbar]}\z^2 
u^{\frac m2-2} |Du|^2\,dx\, d\tau\\ 
&\quad+\frac{\gm}{\sig^2\rho^2}
\iint_{Q\cap[u<\ssigbar]}u^{\frac32m}
\left(\testfunlu\right)^2\,dx\,d\tau,
\end{align*}
where $\gm=4\frac{C_1^2}{C_o}$. Therefore, we conclude that
\begin{align*}
&\frac{C_o}4\iint_{Q\cap[u<\ssigbar]}\z^2 u^{\frac m2-2} 
|Du|^2\,dx\, d\tau\\
&\le\int_{K_{(1+\sig)\rho}\cap[u<\ssigbar]}\z^2(x)
\left(\int_{u(x,0)}^{\ssigbar}\testfunls\,ds\right)dx\\
&\quad+\frac{\gm}{\sig^2\rho^2}\iint_{Q\cap[u<\ssigbar]}
u^{\frac32 m}\left(\testfunlu\right)^2 dx\,d\tau\\
&=J_1+J_2.
\end{align*}
We have
\begin{align*}
&J_1=\int_{K_{(1+\sig)\rho}\cap[u<\ssigbar]}\z^2(x)
\left(\int_{u(x,0)}^{\ssigbar}\testfunls\,ds\right)dx\\
&\le\int_{K_{(1+\sig)\rho}\cap[u<\ssigbar]}
\ssigbar^{1-\frac m2}\left(\int_{u(x,0)}^{\ssigbar}
\frac{\left(\frac{\ssigbar}s\right)^{\frac m2}-1}{m}\,
d\left(\frac s{\ssigbar}\right)\right)dx\\
&=\int_{K_{(1+\sig)\rho}\cap[u<\ssigbar]}
\ssigbar^{1-\frac m2}\left(\int_{\left(\frac{u(x,0)}{\ssigbar}
\right)^m}^1\frac{y^{-\frac12}-1}{m}
\frac{y^{\frac1m-1}}m\,dy\right)dx\\
&=\int_{K_{(1+\sig)\rho}\cap[u<\ssigbar]}
\frac{\ssigbar^{1-\frac m2}}{m^2}
\left[\frac{2m}{2-m}y^{\frac1m-\frac12}-my^{\frac1m}
\right]_{\left(\frac{u(x,0)}{\ssigbar}\right)^m}^{1}dx\\
&=\int_{K_{(1+\sig)\rho}\cap[u<\ssigbar]}
\frac{\ssigbar^{1-\frac m2}}{m}\left[\frac{2}{2-m}
\left(1-\left(\frac{u(x,0)}{\ssigbar}\right)^{1-\frac m2}\right)
-\left(1-\frac{u(x,0)}{\ssigbar}\right)\right]dx\\
&\le\int_{K_{(1+\sig)\rho}\cap[u<\ssigbar]}
\frac{\ssigbar^{1-\frac m2}}{m}
\left[\frac m{2-m}-\frac2{2-m}\frac{u(x,0)}{\ssigbar}
\left(\left(\frac{u(x,0)}{\ssigbar}
\right)^{-\frac m2}-1\right)\right]dx\\
&=\int_{K_{(1+\sig)\rho}\cap[u<\ssigbar]}
{\ssigbar^{1-\frac m2}}\left[\frac1{2-m}
+\frac2{2-m}\frac{u(x,0)}{\ssigbar}
\left(\frac{1-\left(\frac{u(x,0)}{\ssigbar}
\right)^{-\frac m2}}m\right)\right]dx\\
&\le\frac1{2-m}\int_{K_{(1+\sig)\rho}\cap[u<\ssigbar]}
{\ssigbar^{1-\frac m2}}dx\\
&+\frac2{2-m}\int_{K_{(1+\sig)\rho}\cap[u<\ssigbar]}
u(x,0)\frac{u(x,0)^{-\frac m2}-\ssigbar^{-\frac m2}}{m}dx
=J_1^{\prime}+J_1^{\db}.
\end{align*}
\begin{equation*}
J_1^{\prime}=\frac1{2-m}\ssigbar^{1-\frac m2}
\int_{K_{(1+\sig)\rho}\cap[u<\ssigbar]}dx
\le\frac{\gm}{2-m}\ssig^{1-\frac m2}\rho^{N\frac m2},
\quad\text{ where }\ \ \gm=2^N.
\end{equation*}
\begin{align*}
J_1^{\db}&=\frac2{2-m}\int_{K_{(1+\sig)\rho}\cap[u<\ssigbar]}
u(x,0)\frac{u(x,0)^{-\frac m2}-\ssigbar^{-\frac m2}}{m}dx\\
&=\frac2{2-m}\int_{K_{(1+\sig)\rho}\cap[u<\ssigbar]}
\ssigbar^{1-\frac m2}\frac{u(x,0)^{1-\frac m2}}{\ssigbar^{1-\frac m2}}
\frac{\ssigbar^{\frac m2}-u(x,0)^{\frac m2}}{m\ssigbar^{\frac m2}}dx\\
&\le\frac2{2-m}\ssigbar^{1-\frac m2}
\int_{K_{(1+\sig)\rho}\cap[u<\ssigbar]}
\frac{M^{\frac m2}-u(x,0)^{\frac m2}}{m M^{\frac m2}}dx\\
&\le\frac{2\gm}{2-m}\ssigbar^{1-\frac m2}
\Lm_{\frac m2,1}\rho^N=\frac{2\gm}{2-m}\Lm_{\frac m2,1}
\ssig^{1-\frac m2}\rho^{N\frac m2},\quad\text{ where }\ \ \gm=2^N.
\end{align*}
Therefore,
\begin{equation*}
J_1\le\gm(1+\Lm_{\frac m2,1})\rho^{N\frac m2}\ssig^{1-\frac m2}.
\end{equation*}
Moreover,
\begin{align*}
J_2&=\frac{\gm}{\sig^2\rho^2}\iint_{Q\cap[u<\ssigbar]}
u^{\frac32 m}\left(\testfunlu\right)^2 dx\,d\tau\\
&= \frac{\gm}{\sig^2\rho^2}\iint_{Q\cap[u<\ssigbar]}
u^{\frac32 m}\left(\frac{\ssigbar^{\frac m2}-
u^{\frac m2}}{mu^{\frac m2}\ssigbar^{\frac m2}}\right)^2 dx\,d\tau\\
&\le \frac{\gm}{\sig^2\rho^2}\iint_{Q\cap[u<\ssigbar]}
u^{\frac m2}\left(\frac{M^{\frac m2}
-u^{\frac m2}}{mM^{\frac m2}}\right)^2 dx\,d\tau\\
&\le\frac{\gm}{\sig^2\rho^2}\ssigbar^{\frac m2}
\iint_{Q\cap[u<\ssigbar]}\left(\frac{M^{\frac m2}-
u^{\frac m2}}{mM^{\frac m2}}\right)^2 dx\,d\tau\\
&\le\frac{\gm}{\sig^2\rho^2}\ssigbar^{\frac m2}t
\rho^N\sup_{0<\tau<t}\bint_{K_{(1+\sig)\rho}}
\left(\frac{M^{\frac m2}-
u^{\frac m2}}{mM^{\frac m2}}\right)^2 dx\\
&\le\frac{\gm}{\sig^2\rho^2}\Lm_{\frac m2,2}^2
\ssig^{\frac m2}\rho^{N(1-\frac m2)}\,t.
\end{align*}
Hence, we have
\begin{equation}\label{Eq:PM:1}
\begin{aligned}
&\frac{C_o}4\iint_{Q\cap[u<\ssigbar]}\z^2
u^{\frac m2-2}|Du|^2 dx\, d\tau\\
&\le\gm(\Lm_{\frac m2,1}+1)\rho^{N\frac m2}
\ssig^{1-\frac m2}+\frac{\gm}{\sig^2\rho^2}
\Lm_{\frac m2,2}^2\ssig^{\frac m2}\rho^{N(1-\frac m2)}\,t.
\end{aligned}
\end{equation}
Now, if we take the test function 
\begin{equation*}
\vp=\left(\testfunmu\right)_+ \z^2
\end{equation*}
in the weak formulation of (\ref{Eq:4:5})--(\ref{Eq:4:6}) 
and integrate over $Q=K_{(1+\sig)\rho}\times(0,t]$, 
we obtain
\begin{align*}
0=&\iint_{Q}\frac{\pl}{\pl\tau} u 
\left(\testfunmu\right)_+\z^2dx\,d\tau\\
&+\iint_{Q} \bl{A}(x,\tau,u,Du)\cdot 
D\left[\left(\testfunmu\right)_+\z^2\right] dx\,d\tau\\
&= I_3+I_4.
\end{align*}
We estimate these two terms separately. 
\begin{equation*}
\begin{aligned}
I_3&=\iint_{Q} \frac{\pl}{\pl\tau} u 
\left(\testfunmu\right)_+\z^2dx\,d\tau\\
&=\iint_{Q\cap[u>\ssigbar]} \frac{\pl}{\pl\tau} u 
\left(\testfunmu\right)\z^2dx\,d\tau\\
&=\int_{K_{(1+\sig)\rho}\cap[u>\ssigbar]}\z^2(x)
\left(\int_{\ssigbar}^{u(x,t)}\testfunms\,ds\right)dx\\
&-\int_{K_{(1+\sig)\rho}\cap[u>\ssigbar]}\z^2(x)
\left(\int_{\ssigbar}^{u(x,0)}\testfunms\,ds\right)dx.
\end{aligned}
\end{equation*}
Next,
\begin{align*} 
I_4&=\iint_{Q} \bl{A}(x,\tau,u,Du)\cdot 
D\left[\left(\testfunmu\right)_+\z^2\right] dx\,d\tau\\
&=\iint_{Q\cap[u>\ssigbar]} \bl{A}(x,\tau,u,Du)\cdot 
D\left[\left(\testfunmu\right)\z^2\right] dx\,d\tau\\
&=\frac12\iint_{Q\cap[u>\ssigbar]}\z^2 u^{-\frac m2-1}
\bl{A}(x,\tau,u,Du)\cdot Du\,dx\, d\tau\\
&\quad +2\iint_{Q\cap[u>\ssigbar]}\z
\left(\testfunmu\right)\bl{A}(x,\tau,u,Du)\cdot{D\z} dx\, d\tau\\
&\ge \frac{C_o}2\iint_{Q\cap[u>\ssigbar]}\z^2 
u^{-\frac m2-1}u^{m-1}|Du|^2 dx\,d\tau\\
&\quad -2C_1\iint_{Q\cap[u>\ssigbar]}\z
\left(\testfunmu\right) u^{m-1}|Du| |D\z| dx\, d\tau\\
&\ge\frac{C_o}4\iint_{Q\cap[u>\ssigbar]}
\z^2 u^{\frac m2-2} |Du|^2 dx\, d\tau\\
&\quad -\frac{\gm}{\sig^2\rho^2}\iint_{Q\cap[u>\ssigbar]}
u^{\frac32 m}\left(\testfunmu\right)^2 dx\,d\tau,
\end{align*}
where again $\gm=4\frac{C_1^2}{C_o}$. Therefore, we conclude that
\begin{align*}
&\frac{C_o}4\iint_{Q\cap[u>\ssigbar]}\z^2 
u^{\frac m2-2}|Du|^2dx\, d\tau\\
&\le\int_{K_{(1+\sig)\rho}\cap[u>\ssigbar]}
\z^2(x)\left(\int_{\ssigbar}^{u(x,0)}\testfunms\,ds\right)dx\\
&\quad +\frac{\gm}{\sig^2\rho^2}\iint_{Q\cap[u>\ssigbar]}
u^{\frac32 m}\left(\frac{u^{\frac m2}-
\ssigbar^{\frac m2}}{m u^{\frac m2}
\ssigbar^{\frac m2}}\right)^2 dx\,d\tau=J_3+J_4.
\end{align*}
We have
\begin{align*}
&J_3=\int_{K_{(1+\sig)\rho}\cap[u>\ssigbar]}\z^2(x)
\left(\int_{\ssigbar}^{u(x,0)}\testfunms\,ds\right)dx\\
&\le\int_{K_{(1+\sig)\rho}\cap[u>\ssigbar]}
\ssigbar^{1-\frac m2}\left(\int_{\ssigbar}^{u(x,0)}
\frac{1-\left(\frac{\ssigbar}s\right)^{\frac m2}}{m}\,
d\left(\frac s{\ssigbar}\right)\right)dx\\
&=\int_{K_{(1+\sig)\rho}\cap[u>\ssigbar]}
\ssigbar^{1-\frac m2}\left(\int_1^{\left(\frac{u(x,0)}{\ssigbar}
\right)^m}\frac{1-y^{-\frac12}}{m}\frac{y^{\frac1m-1}}m\,dy\right)dx\\
&=\int_{K_{(1+\sig)\rho}\cap[u>\ssigbar]}
\frac{\ssigbar^{1-\frac m2}}{m^2}
\left[my^{\frac1m}-\frac{2m}{2-m}y^{\frac1m-\frac12}
\right]_1^{\left(\frac{u(x,0)}{\ssigbar}\right)^m}dx\\
&=\int_{K_{(1+\sig)\rho}\cap[u>\ssigbar]}
\frac{\ssigbar^{1-\frac m2}}{m}
\left[\left(\frac{u(x,0)}{\ssigbar}-1\right)
-\frac{2}{2-m}\left(\left(\frac{u(x,0)}{\ssigbar}
\right)^{1-\frac m2}-1\right)\right]dx\\
&\le\int_{K_{(1+\sig)\rho}\cap[u>\ssigbar]}
\frac{\ssigbar^{1-\frac m2}}{m}\left[\frac m{2-m}
+\frac2{2-m}\frac{u(x,0)}{\ssigbar}
\left(1-\left(\frac{u(x,0)}{\ssigbar}
\right)^{-\frac m2}\right)\right]dx\\
&=\int_{K_{(1+\sig)\rho}\cap[u>\ssigbar]}{\ssigbar^{1-\frac m2}}
\left[\frac1{2-m}+\frac2{2-m}\frac{u(x,0)}{\ssigbar}
\left(\frac{1-\left(\frac{u(x,0)}{\ssigbar}
\right)^{-\frac m2}}m\right)\right]dx\\
&\le\frac1{2-m}\int_{K_{(1+\sig)\rho}\cap[u>
\ssigbar]}{\ssigbar^{1-\frac m2}}dx\\
&\quad+\frac2{2-m}\int_{K_{(1+\sig)\rho}\cap[u>\ssigbar]}
u(x,0)\frac{\ssigbar^{-\frac m2}-
u(x,0)^{-\frac m2}}{m}dx=J_3^\prime+J_3^{\db}.
\end{align*}
\begin{equation*}
J_3^\prime=\frac1{2-m}\ssigbar^{1-\frac m2}
\int_{K_{(1+\sig)\rho}\cap[u>\ssigbar]}dx
\le\frac{\gm}{2-m}\ssig^{1-\frac m2}
\rho^{N\frac m2}\quad\text{ where }\ \ \gm=2^N.
\end{equation*}
\begin{align*}
J_3^{\db}&=\frac2{2-m}\int_{K_{(1+\sig)\rho}
\cap[u>\ssigbar]}u(x,0)\frac{\ssigbar^{-\frac m2}
-u(x,0)^{-\frac m2}}{m}dx\\
&\le\frac2{2-m}\frac{\ssigbar^{-\frac m2}
-M^{-\frac m2}}{m}\ssig=\frac2{2-m}
\frac{M^{\frac m2}-\ssigbar^{\frac m2}}{m M^{\frac m2}}
\rho^{N\frac m2}\ssig^{1-\frac m2}.
\end{align*}
As in the interval $(0,M]$ the function 
\begin{equation*}
f(s)=\frac{M^{\frac m2}-s^{\frac m2}}{m M^{\frac m2}}
\end{equation*}
is convex, we can apply Jensen's inequality and conclude that
\begin{align*}
&J_3^{\db}\le\frac2{2-m}\rho^{N\frac m2}
\ssig^{1-\frac m2}\bint_{K_{(1+\sig)\rho}}
\frac{M^{\frac m2}-u(x,s_1)^{\frac m2}}{m M^{\frac m2}}dx\\
&\le\gm\Lm_{\frac m2,1}\rho^{N\frac m2}\ssig^{1-\frac m2}.
\end{align*}
Therefore,
\begin{equation*}
J_3\le\gm(1+\Lm_{\frac m2,1})\rho^{N\frac m2}\ssig^{1-\frac m2}.
\end{equation*}
As for $J_4$ we have
\begin{align*}
J_4&=\frac{\gm}{\sig^2\rho^2}\iint_{Q\cap[u>\ssigbar]}
\frac{u^{\frac32 m}}{\ssigbar^{m}}
\left(\frac{u^{\frac m2}-\ssigbar^{\frac m2}}{m u^{\frac m2}}
\right)^2 dx\,d\tau\\
&\le\frac{\gm}{\sig^2\rho^2}\frac1{\ssigbar^{m}}
\iint_{Q\cap[u>\ssigbar]}u^{\frac32 m}
\left(\frac{M^{\frac m2}
-\ssigbar^{\frac m2}}{m M^{\frac m2}}\right)^2 dx\,d\tau\\
&\le\frac{\gm}{\sig^2\rho^2}\frac1{\ssigbar^{m}}
\left(\frac{M^{\frac m2}
-\ssigbar^{\frac m2}}{m M^{\frac m2}}\right)^2\,
t\,\sup_{0<\tau<t}\int_{K_{(1+\sig)\rho}}
u^{\frac32m}(x,\tau)\,dx\\
&\le\frac{\gm}{\sig^2\rho^2}\frac{\rho^{Nm}}{\ssig^{m}}
\left(\frac{M^{\frac m2}
-\ssigbar^{\frac m2}}{m M^{\frac m2}}\right)^2\,t\,
\left(\sup_{0<\tau<t}\int_{K_{(1+\sig)\rho}}
u(x,\tau)\,dx\right)^{\frac32m}\rho^{N(1-\frac32m)}\\
&=\frac{\gm}{\sig^2\rho^2}\ssig^{\frac m2}
\left(\frac{M^{\frac m2}
-\ssigbar^{\frac m2}}{m M^{\frac m2}}
\right)^2\,t\,\rho^{N(1-\frac m2)}\\
&\le\frac{\gm}{\sig^2\rho^2}\Lm_{\frac m2,1}^2
\ssig^{\frac m2}\,t\,\rho^{N(1-\frac m2)},
\end{align*}
where we have taken into account Jensen's inequality once more.
Hence, we have
\begin{equation}\label{Eq:PM:2}
\begin{aligned}
&\frac{C_o}4\iint_{Q\cap[u>\ssigbar]}\z^2 
u^{\frac m2-2}|Du|^2 dx\, d\tau\\
&\le\gm(\Lm_{\frac m2,1}+1)\rho^{N\frac m2}
\ssig^{1-\frac m2}+\frac{\gm}{\sig^2\rho^2}
\Lm_{\frac m2,1}^2\ssig^{\frac m2}\,t\,\rho^{N(1-\frac m2)}.
\end{aligned}
\end{equation}
The lemma follows by combining estimates (\ref{Eq:PM:1}) 
and (\ref{Eq:PM:2}).

The use of $\left(\testfunlu\right)_+ \z^2$ as test function can be justified using $u+\eps$ instead of $u$, and then letting $\eps\to0$.\hfill\bbox
\begin{corollary}\label{Cor:PM:1} 
Let $u$ be a non-negative, local, weak solution to 
the singular equations (\ref{Eq:4:5})--(\ref{Eq:4:6}), 
in $E_T$.  There exists a positive constant $\gm$ 
depending only on the data $\data$, such that for 
all cylinders $K_{4\rho}(y)\times[s,t]\subset E_T$, 
and all $\sig\in(0,1)$, 
\begin{equation*}
\begin{aligned}
\frac1{\rho}
\int_s^t\int_{K_\rho(y)}|\bl{A}(x,\tau,u,Du)|dx\,d\tau
&\le\frac\gm{\sig}(\Lm_{\frac m2,1}^2
+\Lm_{\frac m2,2}^2)^{\frac12}
\left(\frac{t-s}{\rho^\lm}\right)\ssig^{m}\\
&\quad +\gm(1+\Lm_{\frac m2,1})^{\frac12}
\left(\frac{t-s}{\rho^\lm}\right)^{\frac12}\ssig^{\frac{m+1}2}
\end{aligned}
\end{equation*}
\end{corollary}
{\it Proof -} Assume $(y,s)=(0,0)$, and 
let $Q=K_\rho\times(0,t]$. By the structure conditions 
of $\bl{A}$ 
\begin{align*}
&\frac1{\rho}
\int_0^t\int_{K_\rho}
|\bl{A}(x,\tau,u,Du)|dx\,d\tau
\le \frac{C_1}{\rho}\iint_Q u^{m-1}|Du|dx\,d\tau\\
&\le\frac{C_1}{\rho}\left(\iint_Q u^{\frac m2-2}|Du|^2 dx\,
d\tau\right)^{\frac12}
\left(\iint_Q u^{\frac32 m}dx\,d\tau\right)^{\frac12}\\
&\le\gm\frac{C_1}{\rho}
\left[(1+\Lm_{\frac m2,1})\rho^{N\frac m2}\,
\mcl{S}_\sig^{1-\frac m2}+\frac{1}{\sig^2\rho^2}
(\Lm^2_{\frac m2,1}+\Lm^2_{\frac m2,2})\,
\mcl{S}_\sig^{\frac m2}\,t\,\rho^{N(1-\frac m2)}
\right]^{\frac12}\\
&\quad\times \left[t\sup_{0<\tau<t}\int_{K_{(1+\sig)\rho}}
u^{\frac32m}(x,\tau)dx\right]^{\frac12}\\
&\le\gm\frac{C_1}{\rho}\left[(1+\Lm_{\frac m2,1})
\rho^{N\frac m2}\,\mcl{S}_\sig^{1-\frac m2}
+\frac{1}{\sig^2\rho^2}(\Lm^2_{\frac m2,1}
+\Lm^2_{\frac m2,2})\,\mcl{S}_\sig^{\frac m2}\,t\,
\rho^{N(1-\frac m2)}\right]^{\frac12}\\
&\quad\times \left[t\,\ssig^{\frac32m}
\rho^{N(1-\frac32m)}\right]^{\frac12}.
\end{align*}
By simple computations, we conclude.\hfill\bbox
\subsection*{Proof of Proposition~\ref{Prop:4:1}}\label{S:PM:2}
We conclude as in the proof of Proposition~\ref{Prop:2:1}, 
relying on Corollary~\ref{Cor:PM:1},
instead of Corollary~\ref{Cor:Weak:1}.\hfill\bbox


\begin{thebibliography}{99}
\bibitem{HP} M.A. Herrero and M. Pierre, The Cauchy 
problem $\dsty u_t=\Delta u^m$ when $0<m<1$, {\it 
Trans. Amer. Math. Soc., {\bf 291}(1), (1985), 145--158}.
\bibitem{DBGL1} E. DiBenedetto, U. Gianazza and N. Liao, 
On the Local Behavior of Non-Negative Solutions to a 
Logarithmically Singular Equation, {\it Discrete 
Continuous Dynamical Systems Ser. B, 17(6), 2012, 1841--1858.}
\bibitem{DBGL2} E. DiBenedetto, U. Gianazza and N. Liao, 
Logarithmically Singular Parabolic Equations as Limits 
of the Porous Medium Equation, {\it Nonlinear Analysis 
Series A: Theory, Methods \& Applications, 75(12), 2012, 4513--4533.}
\bibitem{DBGV-mono} E. DiBenedetto, U. Gianazza and V. Vespri, 
{\it Harnack's Inequality for Degenerate and Singular Parabolic 
Equations}, Springer Monographs in Mathematics, Springer-Verlag, 
New York, 2012.
\bibitem{DBKV} E. DiBenedetto, Y.C. Kwong and V. Vespri, 
Local Space-analyticity of Solutions of Certain Singular 
Parabolic Equations,  {\it Indiana Univ. Math. J., {\bf40}(2), 
(1991), 741--765.}
\bibitem{KN} D. Kinderlehrer and L. Nirenberg, 
Analyticity at the boundary of solutions of 
nonlinear second order parabolic equations, 
{\it Comm. Pure and Appl. math {\bf XXXI}, (1978), 283--338.}
\bibitem{LSU} O.A. Ladyzenskaya, N.A. Solonnikov and
N.N. Ural'tzeva, {\it Linear and Quasilinear Equations 
of Parabolic Type}, Translations of Mathematical Monographs, 23,
American Mathematical Society, Providence, RI, 1967.
\end{thebibliography}
\bye